# A ONE-DIMENSIONAL COAGULATION-FRAGMENTATION PROCESS WITH A DYNAMICAL PHASE TRANSITION

CÉDRIC BERNARDIN AND FABIO LUCIO TONINELLI

ABSTRACT. We introduce a reversible Markovian coagulation-fragmentation process on the set of partitions of $\{1, \ldots, L\}$ into disjoint intervals. Each interval can either split or merge with one of its two neighbors. The invariant measure can be seen as the Gibbs measure for a homogeneous pinning model [10]. Depending on a parameter $\lambda$, the typical configuration can be either dominated by a single big interval (delocalized phase), or be composed of many intervals of order 1 (localized phase), or the interval length can have a power law distribution (critical regime). In the three cases, the time required to approach equilibrium (in total variation) scales very differently with $L$. In the localized phase, when the initial condition is a single interval of size $L$, the equilibration mechanism is due to the propagation of two "fragmentation fronts" which start from the two boundaries and proceed by power-law jumps.

## CONTENTS

## 1. INTRODUCTION

Coagulation-fragmentation phenomena are often modeled by Markov processes where the configuration at a given time is a set of "fragments", each characterized by a positive number, to be interpreted for instance as a length or a mass. The dynamics then consists in fragmentation events (a fragment breaks into, say, two fragments with a conservation of the total length) and coagulation events (say two fragments coalesce and the lengths add up). Most of the mathematical literature (cf. for instance [3] for a recent review) focuses on *mean field* models, where the rate of coagulation of two fragments is a function only of their lengths $\ell_1, \ell_2$ and does not depend on the two fragments being close or far away in some geometric sense. Exceptions are for instance the one-dimensional models studied in [7, 9, 2]:

there, fragments are seen as intervals of the real line and only neighboring ones can coalesce (in these models, however, fragmentation is not allowed).

In the present work, we consider a (Markovian, continuous-time) one-dimensional process where both fragmentation and coagulation occurs. We introduce the model in a discrete setting: on the interval $\{0, 1, \ldots, L\}$ each site contains either 1 or 0 particles (at sites 0 and $L$ a particle is frozen for all times) and fragments are the intervals between two successive particles. Coagulation of two neighboring fragments (resp. fragmentation) is then interpreted as the disappearance (resp. the creation) of a particle. The transition rates are such that the invariant and reversible measure $\pi_L^\lambda$ is proportional to $\lambda^n$, with $\lambda$ a positive parameter and $n$ the number of fragments, times $\prod_{i=1}^n (\ell_i)^{-\rho-1}$, where $\ell_i$ is the length of the $i^{th}$ fragment and $\rho$ is an exponent in $(0, 1)$. It is well known [10] that there exists a critical value $\lambda_c$ such that, as $L$ diverges: i) if $\lambda < \lambda_c$ the typical configuration for $\pi_L^\lambda$ contains a fragment of length $L - O(1)$ and all the others are finite ii) if $\lambda > \lambda_c$ the typical fragment is of length $O(1)$ and the maximal one is of length $O(\log L)$ and finally iii) for $\lambda = \lambda_c$ the fragments have a power-law tail of exponent $\rho$. Here, we study how the mixing time (i.e. the time to get close to equilibrium in total variation) depends on $L$ in the different phases. In particular, we find that the equilibrium phase transition reflects in a dynamical one: the mixing time is essentially of order $L^\rho$ for $\lambda > \lambda_c$ and much smaller (polylogarithmic in $L$) for $\lambda < \lambda_c$.

That equilibrium phase transitions often have a dynamical counterpart is a well known fact. This is the case for instance for the $d(\geq 2)$-dimensional Ising model with free boundary conditions, whose mixing time is of order $\log L$, like for independent spins, at high temperature and exponential in $L^{d-1}$ below the critical temperature [13]. In that case the reason for the dynamical slowdown is that below the critical temperature the phase space breaks into different valleys (pure phases) separated by high energy barriers. In contrast, for our model at $\lambda > \lambda_c$ there is no multiple-valley phenomenon or coexistence of phases. Instead, the reason for the $L^\rho$ behavior is that in the $\lambda > \lambda_c$ phase, starting with the very far-from-equilibrium configuration with just one fragment, equilibration occurs via the propagation of two "fragmentation fronts" which move from the boundaries to the bulk of the system and proceed by power-law jumps (and therefore faster than ballistically). See the discussion of this heuristics in Section 3.1.

The reason why this model shows a non-trivial phenomenology despite its one-dimensional nature can be understood as follows. It is of course possible to interpret the model as a one-dimensional non-conservative interacting particle system (or spin system). However, in our case the creation-destruction rate at $x$ depends not only on the particle configuration in a finite neighborhood of $x$ but on the location of the first particle to the left and to the right of $x$. In other words, jump rates are very non-local. A related fact is that, in the $\lambda > \lambda_c$ phase, equilibrium particle-occupancy correlation functions decay exponentially (cf. Theorem 1) but the equilibrium in the center of a finite box $\Lambda \subset \{1, \ldots, L\}$, conditioned on the configuration in $\{1, \ldots, L\} \setminus \Lambda$, can depend very strongly on the conditioning even if the box is large (cf. Remark 1 below). If this were not the case, the "fragmentation front" heuristics would fail and the mixing time would be of order $\log L$ even in the $\lambda > \lambda_c$ phase like for usual one-dimensional, finite-range, non-conservative particle dynamics.

In the recent works [5, 6] was considered the heat bath dynamics for a $(1 + 1)$-dimensional polymer model, i.e. the path of a one-dimensional simple random walk, interacting with a defect line whose invariant measure, when projected on the set of polymer-line contacts, is also of the form (2.2) with $\rho = 1/2$. In this context, the $\lambda > \lambda_c$ (resp. $\lambda < \lambda_c$) phase is called the localized (resp. delocalized) phase, a terminology that we will use in the rest of the paper even if we forget the underlying polymer model. Despite this superficial similarity, however, the two dynamical problems are very different: in [5, 6] updates correspond to local moves of the polymer configuration, while (in the polymer language) the updates we consider in this work are very non-local and consist in performing the dynamics directly on the set of contacts. In particular, for the polymer local dynamics the mixing time in the delocalized phase turns out to be polynomial in $L$ and larger than in the localized phase, due to a subtle metastability phenomenon [5, 6].

Finally, let us mention that in [17] the question of the dependence of the mixing time on the system size was analyzed for another coagulation-fragmentation process (the "discrete uniform coagulation-fragmentation process") which is however of mean-field type (but results there are more refined than ours).

In the rest of the paper, we abandon the coagulation-fragmentation language since it is more practical to think in terms of particles and holes.

## 2. The model

2.1. **The equilibrium measure.** Let $K(\cdot)$ be the probability measure on the positive integers defined by

$$K(j) = \mathcal{C}_K \, j^{-(1+\rho)}, \quad j = 1, 2, \ldots \tag{2.1}$$

with $0 < \rho < 1$ and $\mathcal{C}_K$ a positive constant such that $\sum_{j \geq 1} K(j) = 1$. We will comment in Remark 3 that the assumption that $K(\cdot)$ is exactly power law can be to some extent relaxed.

For $L \in \mathbb{N}$ and $\lambda > 0$ let $\pi := \pi_L^\lambda$ be the probability measure on

$$\Omega_L = \{\eta = (\eta_0, \eta_1, \ldots, \eta_L) \in \{0, 1\}^{\{0, 1, \ldots, L\}} : \eta_0 = \eta_L = 1\}$$

defined by

$$\pi_L^\lambda(\eta) = \frac{\lambda^{n(\eta)}}{Z_L(\lambda)} \prod_{j=0}^{n(\eta)} K(x_{j+1} - x_j) \tag{2.2}$$

where the configuration $\eta$ is identified with the set $\{x_0 = 0 < x_1 < \ldots < x_{n(\eta)+1} = L\}$ of sites $x$ occupied by a particle (i.e. $\eta_x = 1$). The number $n(\eta)$ is the number of particles strictly located between 0 and $L$ in the configuration $\eta$. For $\Lambda \subset \{0, \ldots, L\}$ and $\eta \in \Omega_L$, we denote by $\eta_\Lambda \in \{0, 1\}^\Lambda$ the trace of $\eta$ on $\Lambda$. For any probability measure $\mu$ on $\Omega_L$, the marginal of $\mu$ on $\{0, 1\}^\Lambda$ is denoted by $\mu|_\Lambda$.

It is possible to express the partition function $Z_L(\lambda) = \sum_{\eta \in \Omega_L} \pi_L^\lambda(\eta)$ in a more compact way as follows. Let $S = \{S_0, S_1, \ldots\}$ be a renewal process on the integers with $S_0 = 0$ and inter-arrival law $K(\cdot)$, i.e. $\mathbf{P}(S_i - S_{i-1} = j) = K(j)$; denote $\mathbf{P}$ its law and $\mathbf{E}$ the corresponding expectation. Then, one has

$$Z_L(\lambda) = \mathbf{E}\left[\lambda^{\sum_{j=1}^{L-1} \mathbf{1}_{\{j \in S\}}} \mathbf{1}_{\{L \in S\}}\right], \tag{2.3}$$

where $\mathbf{1}_{\{j \in S\}}$ is the indicator function of the event $\{\exists\, i > 0 : S_i = j\}$. One therefore recognizes in $Z_L(\lambda)$ the partition function of a homogeneous pinning model with pinning parameter $\log \lambda$ [10, Chap. 2].

The system undergoes a phase transition at $\lambda = \lambda_c := 1$: for $\lambda > 1$ (localized phase) the partition function grows exponentially with $L$ and the set of sites where $\eta_x = 1$ behaves like a renewal sequence with exponential inter-arrival law; for $\lambda < 1$ (delocalized phase) the partition function tends to zero and the set of sites where $\eta_x = 1$ behaves like a transient renewal process (in particular, the number of particles is a geometric random variable). In the next theorem we collect a few well-known equilibrium facts [10, Chap. 2].

**Theorem 1.**

- **Localized phase.** *For $\lambda > 1$ one has $Z_L(\lambda) \sim \exp(LF(\lambda))$ with $F(\lambda) > 0$ (the free energy). Moreover, the law of $\{0 \leq x \leq L : \eta_x = 1\}$ under $\pi_L^\lambda$ is the same as the law of $\{0 \leq x \leq L : x \in \hat{S}\}$ under $\mathbf{P}_\lambda(\cdot \,|\, L \in \hat{S})$, where $\hat{S}$ is a renewal sequence of law $\mathbf{P}_\lambda$ on the integers with $\hat{S}_0 = 0$ and inter-arrival law*

$$K_\lambda(j) := \mathbf{P}_\lambda(\hat{S}_i - \hat{S}_{i-1} = j) = \frac{e^{-F(\lambda)j}K(j)}{\sum_m e^{-F(\lambda)m}K(m)}.$$

  *Correlations under $\pi_L^\lambda$ decay exponentially fast in space: there exists $c(\lambda) > 0$ such that for every $0 \leq a < b < d \leq L$*

$$|\pi_L^\lambda(\eta_b = 1 | \eta_a = \eta_d = 1) - \pi_L^\lambda(\eta_b = 1)| \leq e^{-c(\lambda)\min((b-a),(d-b))} \tag{2.4}$$

  *and*

$$\pi_L^\lambda(\eta_x = 0 \text{ for every } a < x < b) \leq e^{-c(\lambda)(b-a)}. \tag{2.5}$$

- **Delocalized phase.** *For $\lambda < 1$ one has $Z_L(\lambda) \sim c(\lambda)K(L)$ and the law of $\{0 \leq x \leq L : \eta_x = 1\}$ under $\pi_L^\lambda$ is the same as the law of $\{0 \leq x \leq L : x \in \hat{S}\}$ under $\mathbf{P}_\lambda(\cdot \,|\, L \in \hat{S})$, where $\hat{S}$ is a renewal sequence of law $\mathbf{P}_\lambda$ on the integers with $\hat{S}_0 = 0$ and inter-arrival law*

$$K_\lambda(j) := \mathbf{P}_\lambda(\hat{S}_i - \hat{S}_{i-1} = j) = \lambda K(j).$$

  *The renewal $\hat{S}$ is transient, i.e. $\mathbf{P}_\lambda(\hat{S}_1 = \infty) = 1 - \lambda > 0$. Moreover, there exists $c(\lambda) > 0$ such that for every $1 \leq x \leq L - 1$, one has*

$$\pi_L^\lambda(\eta_x = 1) \leq c(\lambda)\left(\frac{L}{x(L-x)}\right)^{1+\rho} \tag{2.6}$$

  *and*

$$\pi_L^\lambda(n(\eta) \geq k) \leq c(\lambda)\exp(-k/c(\lambda)). \tag{2.7}$$

- **Critical point** *For $\lambda = 1$ one has*

$$Z_L(\lambda = 1) = \mathbf{P}(L \in S) \sim \frac{\rho \sin(\pi\rho)}{\pi}L^{\rho-1}.$$

Here and in the following, $c(\lambda)$ denotes some positive constant depending on $\lambda$ which is not the same at each occurrence. Also, considering $S$ as the set of the renewal times we write $L \in S$ for the event $\{\exists\, i \geq 1 : S_i = L\}$.

**Remark 1.** *Based on* (2.4) *one could be tempted to think that in the localized phase the marginal at site $x$ of $\pi_L^\lambda$ conditioned to the configuration $\eta$ outside a box $\{x - \ell, \dots, x + \ell\}$ would depend weakly in the conditioning, for $\ell$ large. This is false. For instance, it is easy to see that, if $\ell \ll \log L$ and we condition on $\eta_y = 0$ for every $y \notin \{L/2 - \ell, L/2 + \ell\}$, the probability that $\eta_{L/2} = 1$ tends to zero with $L \to \infty$, while the unconditioned $\pi_L^\lambda(\eta_{L/2})$ is bounded away from zero.*

## 2.2. The dynamics.

The integral of a function $f$ with respect to a probability measure $\mu$ is denoted by $\mu(f)$, the covariance between the functions $f$ and $g$ by $\mu(f ; g)$. The continuous-time dynamics $\{\eta(t)\}_{t \geq 0}$ we consider (heat bath dynamics or Gibbs sampler) makes updates which consist in adding or deleting one particle at the time. Its generator is

$$\mathcal{L}f = \sum_{x=1}^{L-1} [\mathcal{Q}_x f - f], \quad f : \Omega_L \to \mathbb{R},$$

where $(\mathcal{Q}_x f)(\eta) = \pi_L^\lambda(f \mid \eta_y, y \neq x)$. In words, the dynamics is described as follows. Each site $1, \dots, L-1$ is equipped with a clock whose rings form i.i.d. Poisson Point Processes of intensity 1. If the clock labeled $x$ rings at some time $t$, we replace $\eta_x(t)$ by a new value sampled from the equilibrium distribution conditioned on the instantaneous configuration outside $x$ at time $t$. It is a standard fact that $\pi_L^\lambda$ is reversible for $\mathcal{L}$.

**Remark 2.** *We consider heat-bath transition rates because we will later need to apply the so-called Peres-Winkler censoring inequalities [15], which have been proved only in this case. These inequalities should presumably hold for more general transition rates, provided that the resulting dynamics is reversible and monotone in the sense of Section 4; the results of the present work would then still hold.*

If we denote by $\eta^x$ the configuration $\eta$ where site $x$ has been flipped,

$$(\eta^x)_z = (1 - \mathbf{1}_{x=z})\eta_z + \mathbf{1}_{x=z}(1 - \eta_x), \tag{2.8}$$

the generator can also be rewritten in the more explicit form

$$(\mathcal{L}f)(\eta) = \sum_{x=1}^{L-1} [(1 - \eta_x)c_x(\eta) + \eta_x d_x(\eta)] [f(\eta^x) - f(\eta)] \tag{2.9}$$

where the "creation rate" $c_x(\eta)$ and the "destruction rate" $d_x(\eta)$ are given by

$$c_x(\eta) = \frac{\lambda}{\lambda + \dfrac{K(x_{i+1} - x_i)}{K(x_{i+1} - x)K(x - x_i)}}, \quad \text{if} \quad x_i < x < x_{i+1}, \tag{2.10}$$

and

$$d_x(\eta) = \frac{1}{1 + \lambda \dfrac{K(x_{i+1} - x_i)K(x_i - x_{i-1})}{K(x_{i+1} - x_{i-1})}}, \quad \text{if} \quad x = x_i. \tag{2.11}$$

We recall that in the previous formulas $x_i, i = 0, \dots, N(\eta) + 1$, denote the positions of the particles in the configuration $\eta$. More generally, if $\Lambda$ is a subset of $\{1, \dots, L-$ 1}, the heat bath dynamics with only sites in $\Lambda$ updated and boundary condition $\tau \in \Omega_L$ outside $\Lambda$ is defined by the generator

$$\mathcal{L}_\Lambda^\tau f = \sum_{x \in \Lambda} [\mathcal{Q}_x f - f], \quad f : \Omega_\Lambda^\tau \to \mathbb{R}, \tag{2.12}$$

where $\Omega_\Lambda^\tau = \{\eta \in \Omega_L : \eta_{\Lambda^c} = \tau_{\Lambda^c}\}$. Here and in the following $\Lambda^c$ stands for the complementary set of $\Lambda$ in $\{1, \ldots, L-1\}$. Of course the probability measure $\pi_L^\lambda(\cdot \mid \eta_{\Lambda^c} = \tau_{\Lambda^c})$ is reversible for $\mathcal{L}_\Lambda^\tau$.

2.3. **Relaxation and mixing time.** The semigroup generated by the Markov process $\{\eta(t)\}_{t \geq 0}$ with generator $\mathcal{L}$, see (2.9), is denoted by $e^{t\mathcal{L}}$ and its Dirichlet form is given by

$$\mathcal{E}_L(f; f) = -\pi_L^\lambda(f\mathcal{L}f) = \frac{1}{2} \sum_{x=1}^{L-1} \sum_{\eta \in \Omega_L} \left((1 - \eta_x)c_x(\eta) + \eta_x d_x(\eta)\right) \left[f(\eta^x) - f(\eta)\right]^2 \pi_L^\lambda(\eta).$$

The corresponding spectral gap is defined by

$$\mathrm{gap} = \inf_{f : \Omega_L \to \mathbb{R}} \frac{\mathcal{E}_L(f; f)}{\pi_L^\lambda(f; f)} \tag{2.13}$$

where the infimum is taken over non-constant functions. The relaxation time $\mathrm{T_{rel}}$ is defined as the inverse of the gap and, for any $f : \Omega_L \to \mathbb{R}$ and $t \geq 0$, we have

$$\pi_L^\lambda(e^{t\mathcal{L}}f; e^{t\mathcal{L}}f) \leq e^{-2t/\mathrm{T_{rel}}} \pi_L^\lambda(f; f).$$

Hence $\mathrm{T_{rel}}$ measures the speed of convergence to equilibrium in $\mathbb{L}^2(\pi_L^\lambda)$-norm. Another natural and widely used way to measure this convergence is with respect to the total variation distance. If $\mu, \nu$ are two probability measures on a finite probability space $E$ the total variation distance between $\mu$ and $\nu$ is defined by

$$\|\mu - \nu\| := \frac{1}{2} \sum_{x \in E} |\mu(x) - \nu(x)| = \sup_{A \subset E} |\mu(A) - \nu(A)| = \inf_{X \stackrel{\mathcal{L}}{\sim} \mu, Y \stackrel{\mathcal{L}}{\sim} \nu} \mathbb{P}(X \neq Y), \tag{2.14}$$

where the infimum is taken over all couplings of $\mu, \nu$ with $X \stackrel{\mathcal{L}}{\sim} \mu$ meaning that the law of the random variable $X$ is $\mu$. The mixing time $\mathrm{T_{mix}}$ is defined by

$$\mathrm{T_{mix}} = \inf\{t \geq 0 : \sup_{\sigma \in \Omega_L} \|\mu_t^\sigma - \pi_L^\lambda\| \leq (2e)^{-1}\}.$$

where $\mu_t^\sigma$ stands for the law at time $t$ of the process starting from the initial configuration $\sigma$. The choice of the numerical factor $(2e)^{-1}$ is irrelevant (any other value smaller than $1/2$ would be essentially equivalent) but with this definition we have

$$\sup_{\sigma \in \Omega_L} \|\mu_t^\sigma - \pi_L^\lambda\| \leq e^{-\lfloor t/\mathrm{T_{mix}} \rfloor}. \tag{2.15}$$

With respect to the spectral gap, the mixing time is much more sensitive to atypical initial configurations (w.r.t. equilibrium). We have the following general bounds between $\mathrm{T_{rel}}$ and $\mathrm{T_{mix}}$

$$\mathrm{T_{rel}} \leq \mathrm{T_{mix}} \leq \log\left(\frac{2e}{\pi^*}\right) \mathrm{T_{rel}} \tag{2.16}$$

where $\pi^* = \min_{\eta \in \Omega_L} \pi_L^\lambda(\eta)$.

## 3. Results

Our first theorem gives the correct order of the mixing time (up to logarithmic corrections) in the localized phase $\lambda > 1$.

**Theorem 2.** *Let $\lambda > 1$. There exist positive constants $C_1(\lambda)$, $C_2(\lambda)$ and $C$ such that for $L \geq 2$*

$$C_1(\lambda)L^\rho \leq T_{\mathrm{mix}} \leq C_2(\lambda)L^\rho(\log L)^{\mathrm{C}}.$$

*The constant $C$ depends only on the exponent $\rho$ in the definition of $K(\cdot)$.*

That the equilibrium phase transition at $\lambda = 1$ has a dynamical counterpart is implied by our second result, which shows that in the delocalized phase the mixing time grows (poly)-logarithmically in $L$.

**Theorem 3.** *Let $\varepsilon > 0$. There exists $C_3 < \infty$ and, for every $\lambda < 1$, positive constants $L_0(\lambda, \varepsilon), C_4(\lambda)$ such that for every $L \geq L_0(\lambda, \varepsilon)$,*

$$(1 - \varepsilon) \log L \leq T_{\mathrm{mix}} \leq (\log L)^{C_3} \tag{3.1}$$

*and $T_{\mathrm{rel}} \leq C_4(\lambda)$ for every $L$. The constant $C_3$ depends only on the exponent $\rho$ in the definition of $K(\cdot)$.*

We conjecture that the mixing time in the delocalized phase is actually $O(\log L)$, see Remark 4.

In the critical case $\lambda = 1$ we can only prove non-optimal bounds both for the relaxation time and for the mixing time, which however are sufficient to show that the relaxation time has a different scaling than in the delocalized phase: one has

$$C_5\, L^\rho \leq T_{\mathrm{mix}} \leq C_6 L^{2+\rho} \tag{3.2}$$

and

$$C_5\, L^\rho \leq T_{\mathrm{rel}} \leq C_6 L^{1+\rho} \tag{3.3}$$

for some positive constants $C_5, C_6$ (cf. Section 7).

A very interesting open problem is to understand whether $T_{\mathrm{rel}}$ diverges with $L$ in the localized phase.

### 3.1. **A bit of heuristics and comments.** 

The reason for the different scaling of the mixing time according to whether $\lambda \gtrless \lambda_c$ can be understood as follows. Take $\lambda > \lambda_c$ and start the dynamics from the configuration without particles in $\{1, \ldots, L-1\}$ (i.e. there is just one fragment), which should be so to speak as far away as possible from equilibrium. Then, particles (i.e. fragments) start to be created at the two endpoints of the system and the two "fragmentation fronts" tend to invade the whole interval $\{1, \ldots, L\}$. The reason for the front phenomenon is that, even if $\lambda$ is large (which favors the presence of particles), creating a particle inside a very large set of empty sites is very unlikely (see (2.10)) and therefore particles tend to be created close to positions where there is a particle already. We will see that the position of, say, the front moving rightwards can be approximated by a continuous time, positive-jump random walk whose transition rate from $x$ to $x + d$ is approximately $1/d^{1+\rho}$. It is then clear from classical results for sums of heavy-tailed random variables that the two fronts will meet in a time of order $L^\rho$. Once the two fronts meet, the system is essentially at equilibrium. Making this heuristics rigorous is the core of Section 5. One technical difficulty is that the front position is not necessarily increasing (due to coagulation events near the front), so the front itself will be defined in a suitably coarse-grained sense. This is the reason for the appearance of the logarithmic factor (which we believe to be spurious) in Theorem 2.

The situation is very different in the delocalized phase. Again start as far away from equilibrium as possible, which however in this case means starting from the configuration with particles at each site ($L$ fragments). Then, particles start disappearing essentially independently everywhere (and not starting from the boundaries) and reappearing less often. This is clear for very small $\lambda$ since the ratio creation rate/destruction rate is of order $\lambda$, see (2.10)-(2.11). Moreover, the latter ratio decreases very quickly during the equilibrium relaxation because the creation rates become smaller and smaller as new empty sites appear. The mixing time, then, should not be very different from the first time when all the Poisson clocks associated to the sites in $\{1, \ldots, L\}$ have rung, which of course is logarithmic in $L$. In reality, however, proving that mixing occurs quickly is much harder. A natural idea would be to apply simple-minded path coupling [11, Sec. 14.2], trying to prove that in a time of order 1 the dynamics contracts the Hamming distance between configurations. This works well for $\lambda$ very small and leads to $\mathrm{T_{mix}} = O(\log L)$, as the reader can check (cf. Remark 4), but has no chance to work up to $\lambda = 1$. The next idea would be to replace the single-flip dynamics with a block dynamics where one updates the particle configurations in blocks whose size depends on $\lambda$ and becomes large as $\lambda \nearrow 1$. In this case, path coupling works but the problem is then to compare the mixing time of the single-flip dynamics with that of the block dynamics. For non-conservative attractive particle systems with finite-range flip rates, it would be easy (applying the so-called Peres-Winkler censoring inequality, see Section 5) to show that the ratio of the two mixing times depends only on $\lambda$ and not on $L$. In our case, however, due to the non-local nature of the flip rates, this does not work. Therefore, we had to devise a different iterative strategy to prove fast mixing in the delocalized phase, that we think can be of independent interest. The drawback is that we get a sub-optimal upper bound on the mixing time (polylogarithmic instead of $O(\log L)$). As a side remark, comparing spectral gaps for single-flip and block dynamics is instead rather standard and this is the reason why we get the optimal result $\mathrm{T_{rel}} = O(1)$ in Theorem 3.

## 4. MONOTONICITY

In this paper a central role is played by monotonicity properties of the dynamics. On $\Omega_L$ there is a natural partial ordering: we say that $\eta \leq \xi$ if $\eta_x \leq \xi_x$ for every $x = 0, \ldots, L$. Analogously, for $\Lambda \subset \{0, \ldots, L\}$ we write $\eta_\Lambda \leq \xi_\Lambda$ if $\eta_x \leq \xi_x$ for every $x \in \Lambda$. The maximal (filled) configuration $\eta_x = 1, 0 \leq x \leq L$, is denoted by "$+$" and the minimal (empty) configuration $\eta_x = 0, x = 1, \ldots, (L-1)$, by "$-$".

Let $\Lambda \subset \{0, \ldots, L\}$ and $\tau, \tau'$ be two boundary conditions. For $\xi \in \Omega_\Lambda^\tau$ and $t > 0$ let $\eta^{\xi;\tau}(t)$ denote the configuration at time $t$ of the dynamics in $\Lambda$, which evolves with boundary condition $\tau$, started from $\xi$. Let also $\mu_t^{\xi;\tau}$ be the law of $\eta^{\xi;\tau}(t)$.

**Lemma 1.** *Let $\Lambda \subset \{0, \ldots, L\}$, let $\tau, \tau' \in \Omega_L$ and $\xi \in \Omega_\Lambda^\tau, \xi' \in \Omega_\Lambda^{\tau'}$. Assume that $\tau \leq \tau'$ and $\xi \leq \xi'$. Then, it is possible to couple the dynamics $\{\eta^{\xi;\tau}(t)\}_{t \geq 0}$ and $\{\eta^{\xi';\tau'}(t)\}_{t \geq 0}$ in such a way that, almost surely, $\eta^{\xi;\tau}(t) \leq \eta^{\xi';\tau'}(t)$ for every $t \geq 0$.*

*Proof.* Let $\alpha \leq a < x < b \leq \beta$. One observes that

$$\frac{\lambda}{\lambda + \dfrac{K(b-a)}{K(b-x)K(x-a)}} \geq \frac{\lambda}{\lambda + \dfrac{K(\beta - \alpha)}{K(\beta - x)K(x - \alpha)}}$$

and therefore

$$\frac{1}{1 + \lambda \dfrac{K(b-x)K(x-a)}{K(b-a)}} \leq \frac{1}{1 + \lambda \dfrac{K(\beta-x)K(x-\alpha)}{K(\beta-\alpha)}}.$$

(With $s = a - \alpha, t = x - a, u = b - x, v = \beta - b$ these inequalities are equivalent to

$$\frac{1}{s+t} + \frac{1}{u+v} \leq \frac{1}{t} + \frac{1}{u}$$

which is trivially satisfied.) Recalling the definition (2.10) for the rate of creation of a particle at $x$, one sees that above inequalities imply that the rate of creation (resp. of destruction) of a particle at an empty (resp. occupied) site $x$ is increasing (resp. decreasing) with respect to $\eta_y, y \neq x$. The claim then easily follows from a standard coupling argument. $\qquad\square$

An immediate consequence is that, under the hypotheses of the Lemma, for every $t$ one has $\mu_t^{\xi;\tau} \preceq \mu_t^{\xi';\tau'}$ where, given two probability measures $\mu, \nu$ on $\Omega_L$, we say that $\mu \preceq \nu$ ($\mu$ is stochastically dominated by $\nu$) if $\mu(f) \leq \nu(f)$ for every function $f$ which is increasing w.r.t. the above specified partial ordering of $\Omega_L$. Letting $t \to \infty$ we also see that if $\tau \leq \tau'$

$$\pi_L^\lambda(\cdot | \eta_{\Lambda^c} = \tau_{\Lambda^c}) \preceq \pi_L^\lambda(\cdot | \eta_{\Lambda^c} = \tau'_{\Lambda^c}). \tag{4.1}$$

It is actually a standard fact that, given $\Lambda$, one can construct *all* the processes $\{\eta^{\xi;\tau}(t)\}_{t \geq 0}$ for different initial conditions $\xi$ and boundary conditions $\tau$ on the same probability space, in such a way that $\eta^{\xi;\tau}(t) \leq \eta^{\xi';\tau'}(t)$ for every $t > 0$, whenever $\tau \leq \tau', \xi \leq \xi'$ (global monotone coupling). In the following, we use this fact implicitly whenever we say "by monotonicity...". The law of the global monotone coupling will be denoted generically as $\mathbb{P}$.

**Remark 3.** *It is not strictly speaking necessary that $K(\cdot)$ is exactly power-law as in (2.1). Indeed, for our results to hold it is enough to require that $K(j) \sim C_K j^{-1-\rho}, \sum_{j \geq 1} K(j) = 1$ and in addition that*

$$\frac{K(b-x)K(x-a)}{K(b-a)} \geq \frac{K(\beta-x)K(x-\alpha)}{K(\beta-\alpha)} \tag{4.2}$$

*whenever $\alpha \leq a < x < b \leq \beta$, which guarantees that Lemma 1 holds. As an example, let $\hat{K}(j)$ be the probability that the first return to 0 of the symmetric simple random walk on $\mathbb{Z}$, started at 0 and conditioned to be non-negative (call $S$ the trajectory of such conditioned random walk and $\mathbf{P}$ its law), occurs at time $2j$. Define then $K(j) = \hat{K}(j)/\sum_{n \geq 1} \hat{K}(n)$ to guarantee that $K$ is normalized to 1. Then, it is known that $K(j) \sim C_K j^{-3/2}$, i.e. $\rho = 1/2$. It is not hard to realize that (4.2) is verified in this case. Indeed, given $\mu > 0$ and positive even integers $a < b$,*

let $\mathbf{P}^\mu_{a,b}$ be the law on $S$ such that

$$\mathbf{P}^\mu_{a,b}(A) = \frac{\mathbf{E}[1_A \, \mu^{\sum_{j=a+1}^{b-1} 1_{S_j=0}} | S_a = S_b = 0]}{\mathbf{E}[\mu^{\sum_{j=a+1}^{b-1} 1_{S_j=0}} | S_a = S_b = 0]}.$$

By the FKG inequalities (for this it is important that the increments of $S$ are $\pm 1$; observe that $\{S_{2x} = 0\}$ is a decreasing event), it follows that

$$\mathbf{P}^\mu_{2a,2b}(S_{2x} = 0) \geq \mathbf{P}^\mu_{2\alpha,2\beta}(S_{2x} = 0). \tag{4.3}$$

Dividing left- and right-hand side by $\mu$ and letting $\mu \searrow 0$, this gives (4.2).

This observation can be generalized to every $\rho \in (0,1)$, using the construction by K. Alexander [1] of (asymmetric) random walks on $\mathbb{Z}^+$, with $\pm 1$ increments, such that the law of their first return to zero behaves like $K(j) \sim C_K j^{-1-\rho}$ (see [1, Th. 2.1]).

## 5. LOCALIZED PHASE

We first give a rough upper bound for the mixing time $\mathrm{T}_{\mathrm{mix}}(\Lambda, x)$ of the heat bath dynamics in an interval $\Lambda \subset \{1, \ldots, L-1\}$ of length $\ell$ when a particle at some site $x \in \Lambda$ is kept alive during the time evolution. Hence we consider the heat bath dynamics with updates in $\Lambda \setminus \{x\}$ and boundary condition $\tau$ such that $\tau_x = 1$. The generator is given by $\mathcal{L}^\tau_{\Lambda \setminus \{x\}}$, see (2.12).

**Lemma 2.** *Let $\alpha = 3 + \rho$. There exists a positive constant $C := C(\lambda)$ independent of $\ell := |\Lambda|$ and $L$ such that*

$$\mathrm{T}_{\mathrm{mix}}(\Lambda, x) \leq C \ell^\alpha.$$

The proof (see Section 5.6) is based on a geometric technique introduced in [18, 8] for bounding the spectral gap of reversible Markov chains, cf. also [16, Ch. 3]. Such technique was applied for instance in [12] to prove that the mixing time of the Glauber dynamics at inverse temperature $\beta$ for the two-dimensional Ising model in a $m \times n$ rectangular box is upper bounded by $m^{c_1} \exp(c_2(\beta)n)$ if $n < m$. The additional difficulty in our model is that the interactions are not finite-range and decay slowly with distance (in contrast with the Ising model case) and that the transition rates are not bounded away from zero.

We start with the proof of the mixing time upper bound in Theorem 2 and postpone the lower bound to Section 5.3. By monotonicity it is easy to show (cf. for instance the proof of Eq. (2.10) in [14]) that

$$\sup_{\sigma \in \Omega_L} \|\mu^\sigma_t - \pi^\lambda_L\| \leq 2L \max\left( \|\mu^+_t - \pi^\lambda_L\|, \|\mu^-_t - \pi^\lambda_L\| \right), \tag{5.1}$$

so we get the desired bound if we show that $\max\left( \|\mu^+_t - \pi^\lambda_L\|, \|\mu^-_t - \pi^\lambda_L\| \right) \leq (4eL)^{-1}$ for some $t = O(L^\rho (\log L)^{\mathrm{C}})$.

Since we are in the localized phase we expect that equilibration occurs faster when starting from the full configuration "$+$" than from the empty one "$-$". Indeed, the next result says that starting from "$+$" the variation distance from equilibrium is smaller than $1/(4eL)$ after time of order $(\log L)^{1+\alpha}$ with $\alpha$ the exponent in Lemma 2 (just choose $\ell = t^{1/(1+\alpha)}$ in the statement of Lemma 3).

**Lemma 3.** *There exist positive constants $C := C(\lambda), c := c(\lambda)$ depending only on $\lambda$ such that, for any $\ell \leq L$,*

$$\|\mu_t^+ - \pi_L^\lambda\| \leq C L \left( e^{-c\ell} + e^{-ct/\ell^\alpha} \right). \tag{5.2}$$

*Proof.* This is rather standard, see for instance the proof of Corollary 1.9 in [14]. Assume for simplicity that $\ell$ is an even integer. If $\{\eta^{eq}(t)\}_{t \geq 0}$ denotes the evolution started from the equilibrium distribution, we have by definition (2.14) of variation distance and then by monotonicity (with $\mathbb{P}$ the law of the global coupling):

$$\|\mu_t^+ - \pi_L^\lambda\| \leq \mathbb{P}(\eta^+(t) \neq \eta^{eq}(t)) \leq \sum_{x=1}^{L-1} \left[ \mathbb{P}(\eta_x^+(t) = 1) - \pi_L^\lambda(\eta_x = 1) \right]. \tag{5.3}$$

For $\ell/2 \leq x \leq L - \ell/2$ we have

$$\mathbb{P}(\eta_x^+(t) = 1) - \pi_L^\lambda(\eta_x = 1) \leq \left[ \mathbb{P}^{\ell,+}(\eta_{\ell/2}(t) = 1) - \pi_\ell^\lambda(\eta_{\ell/2} = 1) \right] + e^{-c\ell}, \tag{5.4}$$

where $\mathbb{P}^{\ell,+}$ denotes the law of the dynamics in $\{1, \ldots, \ell-1\}$ with boundary condition $\eta_0 = \eta_\ell = 1$, started from the maximal configuration. We used the exponential decay of equilibrium correlations (cf. (2.4)) to deduce $|\pi_L^\lambda(\eta_x = 1) - \pi_\ell^\lambda(\eta_{\ell/2} = 1)| \leq \exp(-c\ell)$ and monotonicity of the dynamics to deduce $\mathbb{P}(\eta_x^+(t) = 1) \leq \mathbb{P}^{\ell,+}(\eta_{\ell/2}(t) = 1)$ (recall that $x$ is at distance at least $\ell/2$ from both endpoints of the system). By Lemma 2 (with $\ell = L$) the mixing time $\mathrm{T}_{\mathrm{mix}}(\ell)$ of the dynamics defined by $\mathbb{P}^{\ell,+}$ is $O(\ell^\alpha)$ and we recall (cf. (2.15)) that

$$\left| \mathbb{P}^{\ell,+}(\eta_{\ell/2}(t) = 1) - \pi_\ell^\lambda(\eta_{\ell/2} = 1) \right| \leq e^{-\lfloor t/\mathrm{T}_{\mathrm{mix}}(\ell) \rfloor}.$$

The terms in (5.3) where $x$ is within distance $\ell/2$ from 0 or $L$ are treated similarly. $\qquad\square$

The real job is to estimate the time needed to reach equilibrium starting from the minimal configuration. We turn to this problem now.

### 5.1. The censored dynamics. We start with some definitions. Let

$$\ell := \ell(L) = \lfloor (\log L)^{1+\beta} \rfloor,$$

$\beta$ being a positive fixed parameter and assume for simplicity that $L = (r+1)\ell$ for some $r \geq 2$.

**Definition 1.** *For $j = 0, \ldots, r+1$, we define $y_j = j\ell$ and we note $\Lambda_j = \{y_j - \ell, \ldots, y_j + \ell\} \cap \{1, \ldots, L\}$. For later convenience, we define also "reduced boxes" $\Lambda'_j$ as follows. For $j = 1, \ldots, r$ let*

$$\Lambda'_j := \{x \in \Lambda_j : (1/10)\ell \leq |x - y_j| \leq (9/10)\ell\}$$

*while for $j = 0$ and $j = r+1$ let*

$$\Lambda'_j := \{x \in \Lambda_j : |x - y_j| \leq (9/10)\ell\}.$$

We partition the box $\{0, \ldots, L\}$ into the $r$ overlapping boxes $\Lambda_j, j = 1, \ldots, r$. Note that $y_j$ is the leftmost point of the box $\Lambda_{j+1}$. We also introduce a sequence of deterministic times $T_k = k(1 + \ell^{\alpha+1}), k \geq 0$, with $\alpha$ the exponent in Lemma 2.

Let $\mu'_t$ be the law at time $t$ of the *censored dynamics* $\{\eta'(t)\}_{t \geq 0}$ obtained from the following scheme.

- The initial condition is the minimal configuration "$-$".

- For $t \in [T_k, T_k + 1)$, the dynamics $\{\eta'_t\}_t$ is given by the heat bath dynamics on $\{1, \dots, L-1\}$.
- For $t \in [T_k + 1, T_{k+1})$ the dynamics $\{\eta'_t\}_t$ is obtained by the heat-bath dynamics by keeping only the updates in $\{1, \dots, L-1\} \setminus \{y_j, j = 1, \dots, r\}$.

Observe that $\pi_\lambda^L$ is invariant for the censored dynamics.

**Lemma 4.** *We have for every $t > 0$*

$$\|\mu_t^- - \pi_L^\lambda\| \leq \|\mu'_t - \pi_L^\lambda\|$$

*and $\mu'_t \preceq \mu_t^-$; moreover the functions $\mu_t^-/\pi_L^\lambda$, $\mu'_t/\pi_L^\lambda$ are decreasing.*

This is an immediate consequence of the so-called "censoring inequalities" by Peres and Winkler which say that, for a monotone heat bath dynamics started from the maximal or minimal configuration, censoring (i.e. eliminating) certain updates (according to a schedule which does not depend on the realization of the process) increases the variation distance from equilibrium.

**Theorem 4.** *([15], see also [14, Th. 2.5]) Let $n \in \mathbb{N}$, $0 < t_0 < t_1 < \dots t_n = T$ and $\Lambda_i \subset \{1, \dots, L-1\}$, $i = 1, \dots, n$. Let $\mu_0$ be a law on $\Omega_L$ such that the function $\eta \mapsto \mu_0(\eta)/\pi_L^\lambda(\eta)$ is increasing (resp. decreasing). Let $\mu_T$ be the law at time $T$ of the continuous-time, heat-bath dynamics in $\{1, \dots, L\}$, started from $\mu_0$ at time zero. Also, let $\nu_T$ be the law at time $T$ of the modified dynamics which again starts from $\mu_0$ at time zero, and which is obtained from the above continuous time, heat-bath dynamics by keeping only the updates in $\Lambda_i$ in the time interval $[t_{i-1}, t_i)$ for $i = 1, \dots, n$. Then,*

$$\|\mu_T - \pi_L^\lambda\| \leq \|\nu_T - \pi_L^\lambda\| \tag{5.5}$$

*and $\mu_T \preceq \nu_T$ (resp. $\nu_T \preceq \mu_T$); moreover, the functions $\mu_T/\pi_L^\lambda$ and $\nu_T/\pi_L^\lambda$ are both increasing (resp. decreasing).*

Another very useful consequence of monotonicity is the following lemma, which ensures that, as soon as equilibrium is reached in some box, the system remains at equilibrium forever there. The total variation distance $\|\mu - \nu\|_\Lambda$ in a box $\Lambda \subset \{0, \dots, L\}$ between two probability measures $\mu, \nu$ on $\Omega_L$ is defined as the variation distance between the corresponding marginals on $\Lambda$: $\|\mu - \nu\|_\Lambda := \|\mu|_\Lambda - \nu|_\Lambda\|$.

**Lemma 5.** *Let $\Lambda \subset \{0, \dots, L\}$. The functions $t \mapsto \|\mu_t^+ - \pi_L^\lambda\|_\Lambda$, $t \mapsto \|\mu_t^- - \pi_L^\lambda\|_\Lambda$ and $t \mapsto \|\mu'_t - \pi_L^\lambda\|_\Lambda$ are monotone non-increasing.*

*Proof.* Let us give the proof for instance for $\|\mu'_t - \pi_L^\lambda\|_\Lambda$. For lightness of notation we write here $\pi := \pi_L^\lambda$. For $\eta \in \Omega_L$ write $\eta = (\eta_\Lambda, \eta_{\Lambda^c})$, so that for instance $\pi(\eta_\Lambda) = \pi|_\Lambda(\eta_\Lambda)$. From Theorem 4 we know that the function $\eta \mapsto \mu'_t(\eta)/\pi(\eta)$ is decreasing for every given $t$. We deduce that the same holds for the ratio of the marginals of the two measures:

$$\frac{\mu'_t(\xi_\Lambda)}{\pi(\xi_\Lambda)} \leq \frac{\mu'_t(\xi'_\Lambda)}{\pi(\xi'_\Lambda)} \quad \text{whenever} \quad \xi'_\Lambda \leq \xi_\Lambda. \tag{5.6}$$

Indeed,

$$\begin{aligned}
\frac{\mu'_t(\xi_\Lambda)}{\pi(\xi_\Lambda)} &= \pi\left[\frac{\mu'_t(\xi_\Lambda, \eta_{\Lambda^c})}{\pi(\xi_\Lambda, \eta_{\Lambda^c})} \,\Big|\, \eta_\Lambda = \xi_\Lambda\right] \\
&\leq \pi\left[\frac{\mu'_t(\xi'_\Lambda, \eta_{\Lambda^c})}{\pi(\xi'_\Lambda, \eta_{\Lambda^c})} \,\Big|\, \eta_\Lambda = \xi_\Lambda\right] \leq \pi\left[\frac{\mu'_t(\xi'_\Lambda, \eta_{\Lambda^c})}{\pi(\xi'_\Lambda, \eta_{\Lambda^c})} \,\Big|\, \eta_\Lambda = \xi'_\Lambda\right] = \frac{\mu'_t(\xi'_\Lambda)}{\pi(\xi'_\Lambda)}.
\end{aligned}$$

In the first inequality we used the fact that $\mu'_t(\cdot,\eta_{\Lambda^c})/\pi(\cdot,\eta_{\Lambda^c})$ is decreasing and in the second the fact that $\pi(g(\eta_{\Lambda^c})|\eta_\Lambda = \cdot)$ is decreasing if $g(\cdot)$ is decreasing (this can be seen as a simple consequence of monotonicity of the dynamics in $\Lambda^c$ with respect to the boundary conditions in $\Lambda$). From Theorem 4 it easily follows that, if $s < t$, then $\mu'_s \preceq \mu'_t$ and as a consequence also $\mu'_s|_\Lambda \preceq \mu'_t|_\Lambda$. The claim of Lemma 5 then follows from the following result, applied with $\mu = \mu'_s|_\Lambda, \nu = \mu'_t|_\Lambda$ and $\rho = \pi|_\Lambda$.

**Lemma 6.** *([15, Lemma 2.4]) Let $\rho, \mu, \nu$ be probability measures on a finite, partially ordered probability space. Assume that $\mu \preceq \nu$ and that the function $\eta \mapsto \nu(\eta)/\rho(\eta)$ is decreasing. Then, $\|\nu - \rho\| \leq \|\mu - \rho\|$.*

$\square$

The usefulness of censoring in our case is shown by next result which says that, if $\eta_{y_j} = 1$ at some time $T_k + 1$, then just before time $T_{k+1}$ the system is close to equilibrium in the domain $\Lambda'_j$, cf. Definition 1.

**Proposition 1.** *Fix some $0 \leq j \leq r + 1$ and let $\omega \in \Omega_L$ be a configuration such that $\omega_{y_j} = 1$. There exist positive constants $C, c$ (independent of $j$ and $\ell$) such that*

$$\|\mu'_t(\cdot|\eta'(T_k+1) = \omega) - \pi^\lambda_L\|_{\Lambda'_j} \leq C\,\ell\left(e^{-c(t-T_k-1)/\ell^\alpha} + e^{-c\ell}\right) \tag{5.7}$$

*for every $t \in [T_k + 1, T_{k+1}]$.*

Note that, by definition of $\{T_k\}_k$ and $\ell = \ell(L)$, for $t = T_{k+1}$ the r.h.s. of (5.7) is $O(\exp(-c\ell/2)) = O(L^{-p})$ for every $p > 0$.

*Proof.* We have

$$\|\mu'_t(\cdot \mid \eta'(T_k+1) = \omega) - \pi^\lambda_L\|_{\Lambda'_j} \leq \left\|\mu'_t(\cdot|\eta'(T_k+1) = \omega) - \pi^\lambda_L(\cdot \mid \eta_{y_i} = \omega_{y_i}\forall i)\right\|_{\Lambda'_j}$$
$$+ \left\|\pi^\lambda_L - \pi^\lambda_L(\cdot \mid \eta_{y_i} = \omega_{y_i}\forall i)\right\|_{\Lambda'_j}.$$

From the definition of $\Lambda'_j$, the fact that $\omega_{y_j} = 1$ and the exponential decay of equilibrium correlations (cf. Theorem 1) we have that the second term in the r.h.s. is $O(\exp(-c\ell))$.

The process $\{\eta'(t)\}_t$ in the time interval $[T_k + 1, T_{k+1}]$ is the heat bath dynamics with sites $\{y_i\}_i$ not updated. In particular, conditionally on the event $\{\eta'(T_k+1) = \omega\}$ with $\omega_{y_j} = 1$, the particle on site $y_j$ is alive during the time interval considered. We denote also $\{\tilde\eta^+(t)\}_t$ (resp. $\{\tilde\eta^-(t)\}_t$) the heat bath dynamics with updates only in $\Lambda_j \backslash \{y_j\}$, which starts at time $T_k + 1$ from the configuration which is "+" everywhere in $\{0, \dots, L\}$ (resp. from the configuration which is empty everywhere in $\{1, \dots, L-1\}$ except at $y_j$). Let $\{\eta^\omega(t)\}_{t \geq T_k+1}$ be the evolution which starts at time $T_k + 1$ from the distribution $\pi^\lambda_L(\cdot|\eta_{y_i} = \omega_{y_i}\forall i)$ and where sites $\{y_i\}_i$ are not updated (note that the law of $\eta^\omega(t)$ is $\pi^\lambda_L(\cdot|\eta_{y_i} = \omega_{y_i}\forall i)$ for all $t \geq T_k + 1$). By monotonicity, $\|\mu'_t(\cdot|\eta'(T_k+1) = \omega) - \pi^\lambda_L(\cdot|\eta_{y_i} = \omega_{y_i}\forall i)\|_{\Lambda'_j}$ is upper bounded by the probability that $\eta'(t)$ and $\eta^\omega(t)$ do not coincide in $\Lambda'_j$ (conditionally on $\eta'(T_k + 1) = \omega$). Also, it is obvious (under the global monotone coupling, again) that $(\eta'(t))_{\Lambda'_j} \neq (\eta^\omega(t))_{\Lambda'_j}$ implies $(\tilde\eta^+(t))_{\Lambda'_j} \neq (\tilde\eta^-(t))_{\Lambda'_j}$. Then, thanks to a union bound over the sites of $\Lambda'_j$, one deduces

$$\|\mu'_t(\cdot|\eta'(T_k+1) = \omega) - \pi^\lambda_L(\cdot|\eta_{y_j} = \omega_{y_j}\forall i)\|_{\Lambda'_j} \leq \mathbb{P}((\tilde\eta^+(t))_{\Lambda'_j} \neq (\tilde\eta^-(t))_{\Lambda'_j})$$
$$\leq \sum_{x \in \Lambda'_j} |\tilde\mu^+_t(\eta_x) - \tilde\mu^-_t(\eta_x)|$$

where $\tilde{\mu}_t^{\pm}$ denotes the law of $\tilde{\eta}^{\pm}(t)$.

Note that the equilibrium measure of $\tilde{\eta}^+(t)$ (resp. of $\tilde{\eta}^-(t)$) is $\pi_L^\lambda(\cdot \mid \eta_{\Lambda_j^c \cup \{y_j\}} \equiv 1)$ (resp. $\pi_L^\lambda(\cdot \mid \eta_{\{1,\dots,L-1\} \setminus \Lambda_j} \equiv 0, \eta_{y_j} = 1)$) and by Lemma 2 the mixing time of both dynamics is $O(\ell^\alpha)$. Therefore,

$$
\begin{aligned}
|\tilde{\mu}_t^+(\eta_x) - \tilde{\mu}_t^-(\eta_x)| &\leq C \, e^{-c(t-T_k-1)/\ell^\alpha} \\
&+ \|\pi_L^\lambda(\cdot | \eta_{\Lambda_j^c \cup \{y_j\}} \equiv 1) - \pi_L^\lambda(\cdot | \eta_{\{1,\dots,L-1\} \setminus \Lambda_j} \equiv 0, \eta_{y_j} = 1)\|_{\Lambda_j'}.
\end{aligned}
$$

Thanks to the exponential decay in (2.4), the last term is $O(\exp(-c\ell))$ (here we use that the points in $\Lambda_j^c$ are at distance at least $\ell/10$ from $\Lambda_j'$) and the claim follows. $\qquad\square$

5.2. **The auxiliary process.** We introduce now an auxiliary process $\{\zeta(t)\}_{t \geq 0}$ where $\zeta(t) = (\zeta_1(t), \dots, \zeta_r(t)) \in \{0,1\}^{\{1,\dots,r\}}$ defined from the censored dynamics as follows. The initial condition is $\zeta_j(0) = 0$ for every $j$ and $\zeta(\cdot)$ is constant in the intervals $[T_k, T_{k+1})$. The value of $\zeta_j(T_k)$ is 1 if and only if

- either $\zeta_j(T_{k-1}) = 1$
- or $\zeta_j(T_{k-1}) = 0$ and $\eta'_{y_j}(T_{k-1} + 1) = 1$. This implies, via Proposition 1, that the dynamics $\eta'(t)$ is almost at equilibrium at time $T_k$ in $\Lambda_j'$.

We also set by convention $\zeta_0(t) = \zeta_{r+1}(t) = 1$ for every $t \geq 0$. Observe that $\{\zeta(t)\}_{t \geq 0}$ is a non-Markovian dynamics, but the joint process $\{(\zeta(t), \eta'(t))\}_{t \geq 0}$ is.

Note that in the process $\{\zeta(t)\}_t$ particles never disappear. The next result says that, when there is a particle at $j$ for $\zeta(t)$, with high probability there is a particle for $\eta'(s)$ in $\Lambda_j$, for every $s$ between $t$ and a much larger time, say $L$.

**Lemma 7.** *There exists $c := c(\lambda) > 0$ such that for every $j \in \{0, \dots, r+1\}$ and every $k$ such that $T_k < L$*

$$
\mathbb{P}\left[\zeta_j(T_k) = 1 \text{ and } \exists s \in [T_k, L], \ \eta'(s)|_{\Lambda_j} \equiv 0\right] \leq (1/c) \, e^{-c\ell}.
$$

*Proof.* Let $A_h$ be the event $A_h = \{\eta'_{y_j}(T_{h-1} + 1) = 1\}$ and $\tau = \inf\{h \leq k : A_h \text{ occurs}\}$. We have the following disjoint partition

$$
\{\zeta_j(T_k) = 1\} = \cup_{h \leq k}\{\tau = h\}.
$$

Since $\Lambda_j' \subset \Lambda_j$ and $k < L \ll \exp(c\ell/2)$, the desired claim follows if we show that, for every $h \leq k$,

$$
\mathbb{P}\left[\tau = h \text{ and } \exists s \in [T_k, L], \ \eta'(s)|_{\Lambda_j'} \equiv 0\right] \leq (1/c) \, e^{-c\ell}. \tag{5.8}
$$

The l.h.s. of (5.8) is upper bounded by

$$
\begin{aligned}
&\mathbb{P}\left[\exists s \in [T_k, L] \quad \text{such that } \ \eta'(s)|_{\Lambda_j'} \equiv 0 \ \middle|\ \tau = h\right] \\
&= \mathbb{P}\left[\exists s \in [T_k, L] \quad \text{such that } \ \eta'(s)|_{\Lambda_j'} \equiv 0 \ \middle|\ \eta'_{y_j}(T_{h-1} + 1) = 1\right]
\end{aligned}
$$

where the equality comes from the Markov property of the process $\{\eta'(t)\}_t$. Let

$$
H := \int_{T_k}^L \mathbf{1}(\eta'(s)|_{\Lambda_j'} \equiv 0) \, ds,
$$

so that we need to show that $\mathbb{P}\left[H > 0 \mid \eta'_{y_j}(T_{h-1} + 1) = 1\right] \leq \exp(-c\ell)$. One has via Markov's inequality

$$
\begin{aligned}
\mathbb{P}\left[H > 0 \mid \eta'_{y_j}(T_{h-1} + 1) = 1\right] &\leq \mathbb{P}\left[0 < H < \exp(-a\ell) \mid \eta'_{y_j}(T_{h-1} + 1) = 1\right] \\
&\quad + e^{a\ell}\mathbb{E}\left[H \mid \eta'_{y_j}(T_{h-1} + 1) = 1\right],
\end{aligned}
\tag{5.9}
$$

with $a$ to be chosen. The first term is easily shown to be exponentially small in $\ell$: this is an immediate consequence of the Poisson distribution of clock rings, because $0 < H < \exp(-a\ell)$ requires that there are two clocks that ring within a time interval $\exp(-a\ell)$ and recall that $\exp(-a\ell) = o(L^{-p})$ for every $p$. On the other hand, the second term equals

$$
\int_{T_k}^{L} \mathbb{P}\left(\eta'(s)|_{\Lambda'_j} \equiv 0 \mid \eta'_{y_j}(T_{h-1} + 1) = 1\right) ds.
$$

Thanks to Proposition 1 and Lemma 5 we know that, conditionally on $\eta'_{y_j}(T_{h-1} + 1) = 1$, for every $t > T_k$ the system is within variation distance $\exp(-c_1\ell)$ from equilibrium in the region $\Lambda'_j$, for some $c_1 := c_1(\lambda) > 0$. On the other hand, at equilibrium one has (cf. (2.5))

$$
\pi(\eta'|_{\Lambda'_j} \equiv 0) \leq \exp(-c_2\ell)
$$

for some positive $c_2 := c_2(\lambda)$. In conclusion, the second term in (5.9) is upper bounded by

$$
e^{a\ell}L\left(e^{-c_1\ell} + e^{-c_2\ell}\right)
$$

and the claim follows choosing $a$ smaller than $\min(c_1, c_2)$. $\qquad \square$

It is clear that after some random but finite time (as long as $L < \infty$) one has $\zeta_j = 1$ for every $j \leq r$. The following result quantifies such random time.

**Proposition 2.** *There exist positive constants $c := c(\lambda)$ and $C := C(\rho)$ such that for $R(L) = \lfloor L^\rho \ell^C \rfloor$ we have*

$$
\mathbb{P}(\exists j \in \{1, \ldots, r\},\ \zeta_j(T_{R(L)}) = 0) \leq (1/c)\, e^{-c\ell}. \tag{5.10}
$$

*Proof.* Define the stopping time

$$
\mathcal{T} = \inf\{s > 0 : \exists j \leq (r+1) \text{ such that } \eta'(s)|_{\Lambda'_j} \equiv 0 \text{ and } \zeta_j(s) = 1\}
$$

and the new Markov process $\{(\hat{\zeta}(t), \hat{\eta}'(t))\}_t$ which is constructed like $\{(\zeta(t), \eta'(t))\}_t$ except that all moves, which would lead to $\eta'|_{\Lambda'_j} \equiv 0$ for some $j$ such that at the same time $\zeta_j = 1$, are discarded. Obviously, the two processes can be perfectly coupled until time $\mathcal{T}$. Since from Lemma 7 we see that $\mathbb{P}(\mathcal{T} < L) \leq \exp(-c(\lambda)\ell)$, it is enough to prove (5.10) for the new process.

The process $\{\hat{\zeta}(t)\}_{t\geq 0}$ can also be defined by the following formula

$$
\hat{\zeta}_j(T_k) = \hat{\zeta}_j(T_{k-1}) + (1 - \hat{\zeta}_j(T_{k-1}))\hat{\eta}'_{y_j}(T_{k-1} + 1).
$$

Let $h_j(k) = \max\{n < j : \hat{\zeta}_n(T_k) = 1\}$ denote the position at time $T_k$ of the rightmost particle of $\hat{\zeta}$ before site $j$. Then,

$$
\mathbb{P}(\hat{\eta}'_{y_j}(T_{k-1} + 1) = 1 | (\hat{\zeta}(T_{k-1}), \hat{\eta}'(T_{k-1}))) \geq c_0 \frac{\lambda}{\lambda + 1/K(y_j - y_{h_j(k-1)} + \ell)} \tag{5.11}
$$

for some positive constant $c_0$. Indeed, in the whole time interval $[T_{k-1}, T_{k-1}+1]$ there is a particle in the box $\Lambda'_{h_j(k-1)}$ for the process $\hat{\eta}'$. If $\hat{\eta}'_{y_j}(T_{k-1}) = 1$, then the l.h.s. of (5.11) is lower bounded by the probability that the clock at $y_j$ does not ring before $T_{k-1}+1$. If on the contrary $\hat{\eta}'_{y_j}(T_{k-1}) = 0$, then a lower bound is given by the probability that the clock at $y_j$ rings only once and that a particle is created there. Given the form (2.10) of the creation rate, the fact that $K(\cdot)$ is decreasing and given the definition of $\Lambda'_j$ one obtains the claim (5.11).

Let $\{\tilde{\zeta}(t)\}_{t\geq0}$ be the Markovian process, constant on time intervals $[T_k, T_{k+1})$ where $\tilde{\zeta}(t) = (\tilde{\zeta}_0(t), \tilde{\zeta}_1(t), \dots) \in \{0,1\}^{\{0,1,\dots\}}$, constructed as follows:

- the initial condition is $\tilde{\zeta}_j(0) = 0$, $j \geq 1$ and $\tilde{\zeta}_0(0) = 1$;
- we set $\tilde{\zeta}_0(t) = 1$ for every $t \geq 0$;
- to obtain $\tilde{\zeta}(T_k)$ given $\tilde{\zeta}(T_{k-1})$, sample Bernoulli random variables $B_j(k) \in \{0,1\}, j = 1, 2 \dots$, independent for different $j$, which take the value 1 with probability $\tilde{K}(j - \tilde{h}_j(k-1))$ where

$$\tilde{K}(x) = c_0 \frac{\lambda}{\lambda + 1/K((x+1)\ell)}, \tag{5.12}$$

$c_0$ is the constant appearing in (5.11) and of course $\tilde{h}_j(k) = \max\{n < j : \tilde{\zeta}_n(T_k) = 1\}$. Then,

$$\tilde{\zeta}_j(T_k) = \tilde{\zeta}_j(T_{k-1}) + (1 - \tilde{\zeta}_j(T_{k-1}))B_j(k) \quad j = 1, 2, \dots. \tag{5.13}$$

Observe that $\tilde{K}(\cdot)$ depends on $\lambda$ and $\ell$. Note also that, while $\{\zeta(t)\}_t$ is a particle system on $\{0, \dots, r\}$, $\{\zeta(t)\}_t$ is a particle system on $\{0, 1, \dots\}$. In both cases, particles never disappear once they have appeared.

The reasoning leading to (5.11) shows that $\{\zeta(t)\}_t$ dominates stochastically $\{\tilde{\zeta}(t)\}_t$ so that

$$\mathbb{P}(\exists j \in \{1, \dots, r\}, \ \hat{\zeta}_j(T_{R(L)}) = 0) \leq \mathbb{P}(\exists j \in \{1, \dots, r\}, \ \tilde{\zeta}_j(T_{R(L)}) = 0).$$

In Theorem 5 of next section we will prove that the probability in the r.h.s. is $O(\exp(-c\ell))$, $c := c(\lambda) > 0$, which concludes the proof. □

### 5.3. Front propagation for the process $\{\tilde{\zeta}(t)\}_{t\geq0}$.

In this section we study the process $\tilde{\zeta}(t)$. First we give a result about the propagation of its rightmost particle, which is at position zero at time $t = 0$. Then, we estimate the time it takes for $\tilde{\zeta}(t)$ to fill entirely the interval $\{1, \dots, r\}$ where we recall that $r = L/\ell(L) - 1 = L/\lfloor(\log L)^{1+\beta}\rfloor - 1$.

We associate to the process $\{\tilde{\zeta}(T_k)\}_{k=0,1,\dots}$ the front position, denoted by $f(k)$, defined by

$$f(k) = \sup\{j \geq 0 : \tilde{\zeta}_j(T_k) = 1\}. \tag{5.14}$$

It is easy to see that $\{f(k)\}_{k\geq0}$ is a renewal process on $\{0, 1, \dots\}$ with $f(0) = 0$ and inter-arrival law $\mathbb{P}(f(k) - f(k-1) = j) =: Q(j)$ given by

$$Q(j) = \tilde{K}(j) \prod_{i=j+1}^{\infty} (1 - \tilde{K}(i)), \quad j \geq 1, \quad Q(0) = \prod_{i=1}^{\infty}(1 - \tilde{K}(i)) \tag{5.15}$$

with $\tilde{K}(\cdot)$ defined in (5.12) and we recall that it depends on $\ell = \ell(L)$. Observe that $\sum_{j=0}^{\infty} Q(j) = 1$ (just note that $Q(j) = A(j+1) - A(j)$ if $j > 0$ and $A(j) = \prod_{i=j}^{\infty}(1 - \tilde{K}(i))$).

**Lemma 8.** *There exists $\ell_0 \geq 1$ and $\delta > 0$ such that for $T(n, \ell) = \lfloor n^\rho \ell^{1+\rho} \rfloor$ we have*

$$\forall \ell(L) \geq \ell_0, \forall n \geq 2, \quad \mathbb{P}(f(T(n, \ell)) \in (n/4, 3n/4)) \geq \delta. \tag{5.16}$$

The restriction $n \geq 2$ is just to guarantee that there is at least an integer in $(n/4, 3n/4)$. The proof boils down to showing the convergence in law of $n^{-1} f(T(n, \ell))$ to a non-degenerate stable law of parameter $\rho$ as $n, \ell \to \infty$. This is rather standard and details can be found in Appendix A.

**Theorem 5.** *For any $\varepsilon > 0$ there exists $c(\varepsilon)$ (that can depend also on $\lambda$) such that if $T(L) = L^\rho \ell^{2+\rho+\varepsilon}$ we have for $L$ large*

$$\mathbb{P}(\forall j \in \{1, \ldots, r\}, \; \tilde{\zeta}_j(T(L)) = 1) \geq 1 - e^{-c(\varepsilon)\ell}.$$

*Proof.* Let $\ell \geq \ell_0$ as in Lemma 8. Then, for some positive $\delta$ we have

$$\mathbb{P}\left(\forall j \in \{r/4, \ldots, 3r/4\}, \; \tilde{\zeta}_j(T(r, \ell)) = 0\right) \leq (1 - \delta).$$

Using monotonicity and Markov's property of the process $\tilde{\zeta}$ (monotonicity being guaranteed by the fact that $\tilde{K}(\cdot)$ is decreasing) it follows that at time $t = T(r, \ell)\ell^{1+\varepsilon}$, with probability greater than $1 - (1 - \delta)^{\ell^{1+\varepsilon}}$ we have $\tilde{\zeta}_x(t) = 1$ for some site $x \in \{r/4, \ldots, 3r/4\}$. Conditionally on this event, we can repeat the argument on the two intervals $\{0, \ldots, x\}$ and $\{x, \ldots, r\}$ to claim that we can find two sites $y \in \{x/4, \ldots, 3x/4\}$ and $z \in \{x + (r - x)/4, \ldots, x + 3(r - x)/4\}$ such that $\tilde{\zeta}_{t'}(y) = \tilde{\zeta}_{t'}(z) = 1$ for $t' = T(r, \ell)\ell^{1+\varepsilon} + T(3r/4, \ell)\ell^{1+\varepsilon}$. The probability that all this does occur is bounded below by $1 - (1 - \delta)^{\ell^{1+\varepsilon}} - 2(1 - \delta)^{\ell^{1+\varepsilon}}$.

We repeat this procedure $n$ times, until $r(3/4)^n < 2$ (note that $n = O(\log r) = O(\log L)$). After a time $t = \sum_{j=1}^n T((3/4)^j r, \ell)\ell^{1+\varepsilon}$ and with probability greater than $p = 1 - \sum_{j=1}^n 2^j (1 - \delta)^{\ell^{1+\varepsilon}}$ we have $\tilde{\zeta}_x(t) = 1$ for every $x \in \{0, \ldots, r\}$. The proof is concluded when one observes that $t = O(r^\rho \ell^{2+\rho+\varepsilon}) = O(L^\rho \ell^{2+\varepsilon})$ and $p \geq 1 - L^c(1 - \delta)^{\ell^{1+\varepsilon}} \geq 1 - e^{-c\ell}$ for some positive $c$ which depends on $\varepsilon$ and on $\delta$ (and therefore on $\lambda$). $\qquad\square$

5.4. **Proof of Theorem 2 (upper bound).** In view of (5.1) and Lemma 3, it is sufficient to show

$$\|\mu_t^- - \pi_L^\lambda\| \leq 1/(4eL) \tag{5.17}$$

for some $t = \lfloor L^\rho (\log L)^{C'} \rfloor T_1$ where we recall (see Section 5.1) that

$$T_1 = O((\log L)^{(1+\beta)(\alpha+1)}).$$

By monotonicity, like in (5.3),

$$\|\mu_t^- - \pi_L^\lambda\| \leq \sum_{x=1}^{L-1} \left(\pi_L^\lambda(\eta_x = 1) - \mathbb{P}(\eta_x^-(t) = 1)\right).$$

We consider first the points $x$ in $R := \cup_{j=0}^{r+1} \Lambda_j'$. Thanks to Lemma 4 we have

$$\sum_{x \in R} \left(\pi_L^\lambda(\eta_x = 1) - \mathbb{P}(\eta_x^-(t) = 1)\right) \leq \sum_{x \in R} \left(\pi_L^\lambda(\eta_x = 1) - \mathbb{P}(\eta_x'(t) = 1)\right)$$
$$\leq L \max_j \|\mu_t' - \pi_L^\lambda\|_{\Lambda_j'}. \tag{5.18}$$

From Proposition 2 and the definition $\ell = \lfloor (\log L)^{1+\beta} \rfloor$ we know that, with high probability,
$$\zeta(\lfloor L^\rho (\log L)^{(1+\beta)C} \rfloor) \equiv 1.$$
Via Proposition 1 and Lemma 5 this implies that at time $t$ the dynamics $\eta'(t)$ is within variation distance $\exp(-c\ell)$ from equilibrium in each box $\Lambda'_j$. This shows that the first sum in (5.18) is $o(L^{-p})$ for any $p > 0$.

The contribution from the points in $R' = \{1, \ldots, L-1\} \setminus R$ is treated similarly: one has simply to redefine the points $y_j$ as $y_j = (j+1/2)\ell$ (instead of $y_j = j\ell$) in the definition of the censored dynamics $\{\eta'(t)\}_t$. This way, the union of the new boxes $\cup_{j=1}^{r+1} \Lambda'_j$ covers $R'$ and the argument proceeds similarly as for $R$.

5.5. **Proof of Theorem 2 (lower bound).** We will prove that for times much smaller than $L^\rho$ the dynamics started from the empty configuration is still far from equilibrium (in particular, that there are still very few particles at distance of order $L$ from the boundaries). The idea is to compare the true heat bath dynamics $\{\eta^-(t)\}_{t \geq 0}$ to the dynamics $\{\tilde\eta(t)\}_{t \geq 0}$ (which again starts from the empty configuration) defined as follows:

- the destruction rates $d_x(\eta)$ in (2.11) are set to zero;
- as soon as a particle is created at a site $x$ we put a particle also at the position $L - x$.

Let $\tilde\mu_t$ be the law of $\tilde\eta(t)$. The standard coupling allows us to realize the two processes so that $\tilde\eta(t) \geq \eta(t)$ for every $t$. Since $\|\mu^-_{T_{mix}} - \pi^\lambda_L\| \leq (2e)^{-1}$ and $\pi^\lambda_L\left(\eta|_{\{L/4,\ldots,L/2\}} \equiv 0\right) \to 0$ with $L \to \infty$ (we are in the localized phase), we have for $L$ sufficiently large
$$\mu^-_{T_{mix}}(\eta|_{\{L/4,\ldots,L/2\}} \not\equiv 0) \geq 1/2. \tag{5.19}$$
Since $\tilde\eta(t) \geq \eta(t)$ for every $t$ we have $\tilde\mu_{T_{mix}}(\eta|_{\{L/4,\ldots,L/2\}} \not\equiv 0) \geq 1/2$. We define the front position $\tilde f(t)$ associated to $\tilde\eta(t)$ as
$$\tilde f(t) = \max\{i \in \{0,\ldots,L/2\}, \ \tilde\eta_i(t) = 1\} = L - \min\{i \in \{L/2,\ldots,L\}, \ \tilde\eta_i(t) = 1\}$$
where the equality holds since the process $\tilde\eta(t)$ is by construction symmetric around $L/2$ at all times. Hence for $L$ sufficiently large we have
$$\mathbb{P}\left(L^{-1}\tilde f(T_{mix}) \in [1/4, 1/2]\right) \geq 1/2. \tag{5.20}$$
One can easily prove (see a sketch at the end of this section) that $\tilde f(t)$ is stochastically dominated by a random variable $y(t)$ such that $n^{-1}y(n^\rho)$ converges for $n \to \infty$ to a one sided $\rho$-stable law whose density is denoted $g$. Let $a_L = T_{mix}/L^\rho$ and assume $\liminf_{L \to \infty} a_L = 0$. We have
$$\begin{aligned}
\limsup_L \mathbb{P}\left(L^{-1}\tilde f(T_{mix}) \in [1/4, 1/2]\right) &\leq \limsup_L \mathbb{P}\left(\left(La_L^{1/\rho}\right)^{-1}\tilde f(a_L L^\rho) \geq \frac{a_L^{-1/\rho}}{4}\right) \\
&\leq \limsup_{A \to \infty} \limsup_{n \to \infty} \mathbb{P}\left(\frac{y(n^\rho)}{n} \geq \frac{A}{4}\right) \\
&= \limsup_{A \to \infty} \int_{A/4}^\infty g(x)dx = 0
\end{aligned}$$
which contradicts (5.20). Hence $\liminf a_L > 0$ and the mixing time lower bound is proven.

It remains to give the limit behavior of $\tilde{f}(n^\rho)/n$. Note that $\{\tilde{f}(t)\}_t$ is a Markov process with values in $\{0, \dots, L/2\}$ and positive increments. The rate of the transition from value $\tilde{f}$ to $\tilde{f} + x$ (with $x \leq L/2 - \tilde{f}$) is given (cf. (2.10)) by

$$2\frac{\lambda}{\lambda + K(L - 2\tilde{f})/(K(x)K(L - 2\tilde{f} - x))} \leq \frac{C(\lambda)}{x^{1+\rho}}, \tag{5.21}$$

where the factor 2 is due to the fact that, if a particle appears at $x$ it appears also at $L - x$ for the process $\{\tilde{\eta}(t)\}_t$. Therefore, $\tilde{f}(t) \leq y(t)$ where $\{y(t)\}_t$ is a Markov process on $\{0, 1, \dots\}$ with jump rate $f \to f + x$ given by the r.h.s. of (5.21). Finally, the convergence in law of $\{y(n^\rho)/n\}_n$ to a $\rho$-stable variable is immediate. $\qquad\square$

5.6. **Proof of Lemma 2.** If $x$ is not one of the two endpoints of $\Lambda$, then the dynamics in $\Lambda \cap \{y : y < x\}$ is independent of the dynamics in $\Lambda \cap \{y : y > x\}$. We can therefore assume without loss of generality that $x$ is the left endpoint of $\Lambda$ and moreover by translation invariance we let $x = 0$, so that $\Lambda \setminus \{x\} = \{1, \dots, \ell - 1\}$. Let also $\ell + M = \min\{y \geq \ell : \tau_y = 1\}$.

Recall that $\mathrm{T}_{\mathrm{mix}} \leq \log\left(\frac{2e}{\pi^*}\right) \mathrm{T}_{\mathrm{rel}}$. Here $\pi^*$ is the minimum over $\sigma \in \{0, 1\}^{\{1, \dots, \ell - 1\}}$ of

$$\hat{\pi}^\tau(\sigma) := \pi^\tau_{\Lambda \setminus \{x\}}(\sigma) = \pi^\lambda_L\left(\eta|_{\{1, \dots, \ell - 1\}} = \sigma \,\Big|\, \eta_y = 0, \ell \leq y < \ell + M, \, \eta_{\ell + M} = 1\right) \tag{5.22}$$

It is easy to see that

$$-\log \pi^* \leq C_1(\lambda) \times \ell, \tag{5.23}$$

hence $\mathrm{T}_{\mathrm{mix}} \leq C_2(\lambda) \times \ell \times \mathrm{T}_{\mathrm{rel}}$. To prove (5.23), given $\sigma \in \{0, 1\}^{\{1, \dots, \ell - 1\}}$ let $d$ denote the rightmost particle of $\sigma$ before site $\ell$ (set $d = 0$ if there is no particle). One has then from (5.22) and from (2.2):

$$
\begin{aligned}
\hat{\pi}^\tau(\sigma) &\geq c(\lambda)\frac{K(M + \ell - d)e^{-c'(\lambda)d}}{\sum_{d'=0}^{\ell - 1} K(M + \ell - d')Z_{d'}(\lambda)} \\
&\geq c(\lambda)\frac{K(M + \ell - d)}{K(M + 1)}e^{-c'(\lambda)\ell} \geq c(\lambda)\frac{K(M + \ell)}{K(M + 1)}e^{-c'(\lambda)\ell}
\end{aligned}
$$

since $K(\cdot)$ is decreasing and $Z_L(\lambda)$ grows at most exponentially in $L$, see Theorem 1. The ratio in $K(M + \ell)/K(M + 1)$ is clearly lower bounded by $C \times \ell^{-1-\rho}$ and this concludes the proof of (5.23).

We use a path combinatorics argument to get an upper bound on the relaxation time. For any $\eta, \eta' \in \Omega^\tau_L = \{\eta \in \Omega_L : \eta_{\{1, \dots, \ell - 1\}^c} = \tau_{\{1, \dots, \ell - 1\}^c}\}$ let $x_i, i = 1, \dots, p$ denote the set of ordered sites in $\{1, \dots, \ell - 1\}$ such that $\eta_x \neq \eta'_x$. Consider the path $\gamma(\eta, \eta')$ in the configuration space, of the form $(\eta_0, \dots, \eta_{p-1})$ with $\eta_0 = \eta, \eta_p = \eta'$, $\eta_{i+1} = \eta_i^{x_i}$ (cf. (2.8)). The number $p$ is called the length of the path and is denoted by $|\gamma(\eta, \eta')|$. Plainly, one has $|\gamma(\eta, \eta')| \leq \ell$. Let $\mathcal{A}$ be the set $\{(\eta, \eta^z); \eta \in \Omega^\tau_L, 1 \leq z < \ell\}$ and for $e = (\sigma, \sigma^z) \in \mathcal{A}$ let

$$Q(e) = \hat{\pi}^\tau(\sigma) \times \begin{cases} c_z(\sigma) & \text{if } \sigma_z = 0 \\ d_z(\sigma) & \text{if } \sigma_z = 1. \end{cases} \tag{5.24}$$

By Theorem 3.2.1 of [16] we have

$$T_{\text{rel}} \leq \max_{e \in \mathcal{A}} \left\{ \frac{1}{Q(e)} \sum_{\substack{\eta, \eta' \in \Omega_L^\tau, \\ \gamma(\eta, \eta') \ni e}} |\gamma(\eta, \eta')| \hat{\pi}^\tau(\eta) \hat{\pi}^\tau(\eta') \right\}. \tag{5.25}$$

Consider for definiteness the case $\sigma_z = 1$, the other case being essentially identical by reversibility. One has that $d_z(\sigma)$ is bounded away from zero, uniformly in $\sigma, L, \tau$, so it remains to prove that

$$\frac{1}{\hat{\pi}^\tau(\sigma)} \sum_{\eta, \eta' \in \Omega_L^\tau, \gamma(\eta, \eta') \ni (\sigma, \sigma^z)} \hat{\pi}^\tau(\eta) \hat{\pi}^\tau(\eta') \leq C \ell^{1+\rho}.$$

From the way the path $\gamma(\eta, \eta')$ was constructed, we see that $\gamma(\eta, \eta') \ni (\sigma, \sigma^z)$ implies that $\eta_i$ coincides with $\sigma_i$ for $i \geq z$ while $\eta_i' = \sigma_i$ for $i < z$ and $\eta_z' = 0 = 1 - \sigma_z$. One can identify the configuration $\sigma$ with the positions of its particles up to position $\ell + M$, $0 = x_0^\sigma < \cdots < x_{n(\sigma)+1}^\sigma = \ell + M$ and similarly for $\eta, \eta'$. Call $a$ the index such that $x_a^\sigma = z$, $b$ the index such that $x_b^\eta = z$ and $c$ the largest index such that $x_c^{\eta'} < z$ (observe that $x_c^{\eta'} = x_{a-1}^\sigma$). One then sees that

$$\frac{\hat{\pi}^\tau(\eta) \hat{\pi}^\tau(\eta')}{\hat{\pi}^\tau(\sigma)} = \frac{1}{\hat{Z}^\tau} \frac{\prod_{j=0}^{b-1} [\lambda K(x_{j+1}^\eta - x_j^\eta)] \prod_{j=c}^{n(\eta')} [\lambda K(x_{j+1}^{\eta'} - x_j^{\eta'})]}{\lambda K(z - x_{a-1}^\sigma)} \tag{5.26}$$

$$\leq C \ell^{1+\rho} \frac{\prod_{j=0}^{b-1} [\lambda K(x_{j+1}^\eta - x_j^\eta)] \prod_{j=c}^{n(\eta')} [\lambda K(x_{j+1}^{\eta'} - x_j^{\eta'})]}{\hat{Z}^\tau} \tag{5.27}$$

where

$$\hat{Z}^\tau = \sum_{\sigma \in \Omega_L^\tau} \prod_{j=0}^{n(\sigma)} [\lambda K(x_{j+1}^\eta - x_j^\eta)].$$

Since $K(\cdot)$ is decreasing, the expression in (5.27) increases if we replace $x_c^{\eta'}$ by $z$. Finally, it is immediate to see that

$$\sum_{\substack{\eta, \eta' \in \Omega_L^\tau: \\ \gamma(\eta, \eta') \ni (\sigma, \sigma^z)}} \frac{\prod_{j=0}^{b-1} [\lambda K(x_{j+1}^\eta - x_j^\eta)] [\lambda K(x_{c+1}^{\eta'} - z)] \prod_{j=c+1}^{n(\eta')} [\lambda K(x_{j+1}^{\eta'} - x_j^{\eta'})]}{\hat{Z}^\tau}$$

is equal to $\hat{\pi}^\tau(\eta_z = 1) \leq 1$, which concludes the proof.                    $\square$

## 6. DELOCALIZED PHASE

In this section, we prove Theorem 3. The first subsection is devoted to bounding the mixing time and the second one to bounding the relaxation time.

6.1. **Mixing time.** The lower bound $T_{\text{mix}} \geq (1 - \varepsilon) \log L$ is pretty easy. Let $1 \ll a_L \ll L$ and let $S$ be the first time such that all clocks labeled $a_L \leq i \leq L - a_L$ have rung. Since, by (2.7), $\pi_L^\lambda(\eta_i = 0$ for every $a_L \leq i \leq L - a_L) \to 1$ when $L \to \infty, a_L \to \infty$, one has

$$\|\mu_t^+ - \pi_L^\lambda\| \geq \pi_L^\lambda(\eta_i = 0 \; \forall \; a_L \leq i \leq L - a_L) - \mathbb{P}(\eta_i^+(t) = 0 \; \forall \; a_L \leq i \leq L - a_L)$$
$$\geq 1 + o(1) - \mathbb{P}(S \leq t). \tag{6.1}$$

Note that $S$ has the same law as the maximum of $L - 2a_L$ independent standard exponential variables so that $\mathbb{P}(S \leq (1-\varepsilon)\log L) = (1 - L^{-(1-\varepsilon)})^{L-2a_L}$ which goes to 0 whenever $\varepsilon > 0$ and the mixing time lower bound follows from (6.1).

In order to get the upper bound on $\mathrm{T}_{\mathrm{mix}}$ in Theorem 3, we generalize a bit our model. We define $\Omega_{\mathbb{Z}} := \{0,1\}^{\mathbb{Z}}$ and, given $\tau \in \Omega_{\mathbb{Z}}$, $\Lambda$ an interval in $\mathbb{Z}$, we consider the dynamics in $\Lambda$ with boundary conditions $\tau$ on $\Lambda^c$. If $\tau_x = 0$ for every $x$ the dynamics is defined as the limit as $n \to \infty$ of a dynamics with boundary condition $\tau^{(n)}$ which equals 1 for $|x| > n$ and 0 otherwise (the limit exists by monotonicity of the creation/destruction rates w.r.t. the boundary conditions) and the equilibrium measure is concentrated on the empty configuration. Recall the notation $\mu_t^{\sigma;\tau}$ from Section 4 for the law of the dynamics started from $\sigma$, with boundary condition $\tau$.

The proof of the mixing time upper bound is based on an iterative procedure (on the size of $\Lambda$), together with a path coupling argument by Bubley and Dyer [4] (see also [11, Sec. 14.2]).

Given two configurations $\eta, \eta' \in \Omega_{\Lambda} = \{0,1\}^{\Lambda}$, we introduce the Hamming distance $d(\eta, \eta') = \sum_{i \in \Lambda} |\eta_i - \eta'_i| \leq |\Lambda|$ which counts the number of discrepancies. Given two random laws $\mu, \nu$ on $\Omega_{\Lambda}$, define their Kantorovich distance as

$$d_K(\mu, \nu) := \inf_{\substack{(X,Y): \\ X \overset{\mathcal{L}}{=} \mu, Y \overset{\mathcal{L}}{=} \nu}} \mathbb{E}\, d(X,Y).$$

Then one has the following theorem, which is a direct application of a result from [4]:

**Theorem 6.** *Let the interval $\Lambda$ and the boundary condition $\tau$ be fixed. If there exists $\gamma > 0$ and $t > 0$ such that*

$$d_K(\mu_t^{\sigma;\tau}, \mu_t^{\sigma';\tau}) \leq e^{-\gamma} \tag{6.2}$$

*for every $\sigma, \sigma'$ such that $d(\sigma, \sigma') = 1$, one has*

$$d_K(\mu_{nt}^{\sigma;\tau}, \mu_{nt}^{\sigma';\tau}) \leq e^{-\gamma n} d(\sigma, \sigma') \tag{6.3}$$

*for every $n \in \mathbb{N}$ and every couple $(\sigma, \sigma')$.*

Let us fix some $\beta > 0$ and define

$$T_{\beta}(\Lambda, \tau) := \inf \left\{ t > 0 \; : \; \sup_{d(\sigma,\sigma')=1} d_K(\mu_t^{\sigma;\tau}, \mu_t^{\sigma';\tau}) \leq \frac{1}{|\Lambda|^{\beta}} \right\} \tag{6.4}$$

and

$$\mathbf{T}_{\beta}(\Lambda) := \sup_{\tau} T_{\beta}(\Lambda, \tau) \tag{6.5}$$

(actually $\mathbf{T}_{\beta}(\Lambda)$ depends only on $|\Lambda|$, so we will write rather $\mathbf{T}_{\beta}(|\Lambda|)$).

We will prove the following.

**Proposition 3.** *Fix $1/(1+\rho) < \zeta < 1$ and $0 < \beta < \zeta(1+\rho) - 1$ sufficiently small and let*

$$\ell(L) := \left\lfloor c_2(\lambda) L^{\zeta} \right\rfloor. \tag{6.6}$$

*There exists $C_4(\rho)$ and an integer $L_0(\lambda)$ such that for $L \geq L_0(\lambda)$ one has*

$$\mathbf{T}_{\beta}(L) \leq C_4(\rho) \mathbf{T}_{\beta}(\ell(L)). \tag{6.7}$$

Let us show how this implies the mixing time upper bound

*Proof of Theorem 3 (upper bound).* If $\ell^{(n)}(\cdot)$ denotes the application of the map $x \mapsto \ell(x)$ $n$ times, one has $\ell^{(n)}(L) \le L_0(\lambda)$ for $n \le c(\rho) \log \log L$, for some finite $c(\rho)$, if $L$ is large enough. As a consequence,

$$\mathbf{T}_\beta(L) \le S := (\log L)^{C_5(\rho)} C_6(\lambda)$$

where essentially $C_6(\lambda) = \mathbf{T}_\beta(L_0(\lambda))$.

Since the maximal Hamming distance between configurations in $\Lambda$ is $|\Lambda|$ and thanks to Theorem 6 one has, for the dynamics in $\Lambda$ and irrespective of the boundary condition $\tau$ and on the initial conditions $\sigma, \sigma'$,

$$d_K(\mu_{nS}^{\sigma;\tau}, \mu_{nS}^{\sigma';\tau}) \le \frac{1}{2e} \tag{6.8}$$

if $n > 1/\beta$, where $\beta$ enters the definition (6.4) of $T_\beta(\Lambda, \tau)$. On the other hand, for every $t > 0$ and $\sigma \in \Omega_\Lambda$, denoting for simplicity $\mu_t^{\sigma;\tau}$ by $\mu_t^\sigma$ and by $\pi$ the equilibrium measure in $\Lambda$ with boundary condition $\tau$, one has

$$\|\mu_t^\sigma - \pi\| \le \sup_{\sigma, \sigma'} d_K(\mu_t^\sigma, \mu_t^{\sigma'}). \tag{6.9}$$

Indeed,

$$
\begin{aligned}
\|\mu_t^\sigma - \pi\| &= \frac{1}{2} \sum_\eta \left| \sum_{\eta'} \pi(\eta') \left( \mu_t^\sigma(\eta) - \mu_t^{\eta'}(\eta) \right) \right| \le \frac{1}{2} \sum_{\eta, \eta'} \pi(\eta') \left| \mu_t^\sigma(\eta) - \mu_t^{\eta'}(\eta) \right| \\
&\le \sup_{\eta', \sigma} \|\mu_t^\sigma - \mu_t^{\eta'}\|.
\end{aligned}
$$

By definition of the total variation distance and Markov's inequality, one has

$$\|\mu_t^\sigma - \mu_t^{\eta'}\| \le \mathbb{P}(X \ne Y) = \mathbb{P}(d(X, Y) \ge 1) \le \mathbb{E}(d(X, Y)) \tag{6.10}$$

for every couple of random variables $(X, Y)$ such that $X \overset{\mathcal{L}}{\sim} \mu_t^\sigma, Y \overset{\mathcal{L}}{\sim} \mu_t^{\eta'}$. Choosing $(X, Y)$ as those which realize the infimum in the definition of Kantorovich distance between $\mu_t^\sigma$ and $\mu_t^{\eta'}$, one one obtains (6.9). From (6.8) we deduce that $\|\mu_{nS}^\sigma - \pi\| \le 1/(2e)$ for every initial condition $\sigma$, i.e. the mixing time is upper bounded by $nS$. This proves the upper bound in (3.1) (actually, a stronger statement which is uniform in the boundary conditions). $\qquad \square$

*Proof of Proposition 3.* Consider $\sigma, \sigma' \in \Lambda$ (with $|\Lambda| = L$; for definiteness, say $\Lambda = \{1, \ldots, L\}$) such that $d(\sigma, \sigma') = 1$ and assume for definiteness that $\sigma$ contains one particle more than $\sigma'$, so that in particular $\sigma \ge \sigma'$. Fix some boundary condition $\tau$ outside $\Lambda$ and, for every $\sigma \in \Omega_\Lambda$, denote for simplicity $\mu_t^{\sigma;\tau}$ by $\mu_t^\sigma$. We have by definition of Kantorovich distance and monotonicity

$$d_K(\mu_t^\sigma, \mu_t^{\sigma'}) \le \mathbb{E} d(\eta^\sigma(t), \eta^{\sigma'}(t)) = \sum_{x \in \Lambda} \left( \mu_t^\sigma(\eta_x) - \mu_t^{\sigma'}(\eta_x) \right) \tag{6.11}$$

where $\mathbb{P}$ denotes the global monotone coupling (so that $\eta^\sigma(t) \ge \eta^{\sigma'}(t)$ for every $t$). For simplicity we split the analysis into three cases. Fix some $g$ such that $\zeta > g > \beta/\rho$ and let $-a = \max\{y \le 0 : \tau_y = 1\}$, $b = \min\{y > L : \tau_y = 1\}$.

(i) First case: $\ell(L)/2 < x < L - \ell(L)/2$. Consider the dynamics on $\tilde{\Lambda} := \{x - \ell(L)/2, \ldots, x + \ell(L)/2 - 1\}$ (assume for simplicity that $\ell(L)$ is even and note that $\tilde{\Lambda} \subset \Lambda$) with boundary condition $\tilde{\tau}_y = 1$ for $y \notin \tilde{\Lambda}$. Call $\tilde{\mu}_t^\xi$ its law at time $t$ when the initial condition is $\xi \in \Omega_{\tilde{\Lambda}}$ and call $\tilde{\pi}$ its invariant measure (which is nothing but a space translation of $\pi^\lambda_{\ell(L)}$). One has by monotonicity

$$\mu^\sigma_t(\eta_x) - \mu^{\sigma'}_t(\eta_x) \leq \mu^\sigma_t(\eta_x) \leq \tilde{\mu}^+_t(\eta_x) \leq \tilde{\pi}(\eta_x) + \|\tilde{\mu}^+_t - \tilde{\pi}\|. \tag{6.12}$$

Now choose $t = C_4 \, \mathbf{T}_\beta(\ell(L))$ with $C_4 = C_4(\rho) > 1$. The first term in the r.h.s. of (6.12) is upper bounded by $c(\lambda)\ell(L)^{-1-\rho}$ (cf. (2.6)) and the second by $\ell(L)^{-C_4\beta}$ (thanks to (6.9), the definition of $\mathbf{T}_\beta$ and Theorem 6). These bounds are uniform in the boundary condition $\tau$. As a consequence,

$$\sum_{\ell(L)/2 < x < L - \ell(L)/2} \left( \mu^\sigma_t(\eta_x) - \mu^{\sigma'}_t(\eta_x) \right) < \frac{1}{2L^\beta} \tag{6.13}$$

if $C_4$ is chosen sufficiently large and $L$ is large enough (we used that $\beta < \zeta(1+\rho) - 1$).

(ii) Case $\max(L^g - a, 1) \leq x \leq \ell(L)/2$ (the case $\max(L^g - (b-L), 1) \leq (L-x) \leq \ell(L)/2$ being analogous). This is not much different from Case (i). Note that $\max(L^g - a, 1) \ll \ell(L)/2$ since we assumed $g < \zeta$.

Let this time $\tilde{\Lambda} \subset \Lambda$ be the interval of length $\ell(L)$, whose leftmost border is $\max(-a, x - \ell(L)/2)$ and consider the dynamics in $\tilde{\Lambda}$, with boundary condition $\tilde{\tau}_y = 1$ for $y \notin \tilde{\Lambda}$. Again one has (6.12) where last term in the r.h.s. is negligible if $C_4$ is chosen large enough. As for the first term, one has from (2.6)

$$\tilde{\pi}(\eta_x) \leq \frac{c(\lambda)}{\min((a+x), \ell)^{1+\rho}} \leq \frac{c(\lambda)}{(a+x)^{1+\rho}} + \frac{c(\lambda)}{\ell^{1+\rho}}.$$

Summing on $\max(L^g - a, 1) \leq x \leq \ell(L)/2$ one obtains a quantity which is $O(\ell^{-\rho}) + O(L^{-g\rho}) = o(L^{-\beta})$ due to the way $\beta$ and $g$ were chosen. In particular, for $L$ large enough we have

$$\sum_{\max(L^g - a, 1) \leq x \leq \ell(L)/2} \left( \mu^\sigma_t(\eta_x) - \mu^{\sigma'}_t(\eta_x) \right) < \frac{1}{8L^\beta}. \tag{6.14}$$

(iii) Case $0 < x < \max(L^g - a, 1)$ (or by symmetry $(L-x) \leq \max(L^g - (b-L), 1)$). This situation requires more refined equilibrium estimates. The rough bound $\mu^\sigma_t(\eta_x) - \mu^{\sigma'}_t(\eta_x) \leq \mu^\sigma_t(\eta_x)$ is not sufficient any more and one needs to show that the difference between the two terms is small.

We can clearly assume that $a < L^g$. We consider this time two dynamics on $\hat{\Lambda} = \{1, \ldots, \ell(L)\}$: the first has boundary condition $\tau^{(1)}_x = \tau_x$ for $x \leq 0$ and $\tau^{(1)}_x = 1$ for $x > \ell(L)$ while the second one has boundary condition $\tau^{(2)}_x = \tau_x$ for $x \leq 0$ and $\tau^{(2)}_x = 0$ for $x > \ell(L)$. Call $\mu^{(i),\xi}_t$ $(i = 1, 2)$ their law at time $t$ when the initial condition is $\xi \in \Omega_{\hat{\Lambda}}$ and $\pi^{(i)}, i = 1, 2$ their invariant measures. Then, by monotonicity,

$$\mu^\sigma_t(\eta_x) - \mu^{\sigma'}_t(\eta_x) \leq \mu^{(1),+}_t(\eta_x) - \mu^{(2),-}_t(\eta_x)$$
$$\leq \left[ \pi^{(1)}(\eta_x) - \pi^{(2)}(\eta_x) \right] + \|\mu^{(1),+}_t - \pi^{(1)}\| + \|\mu^{(2),+}_t - \pi^{(2)}\|.$$

Once more, the last two terms are smaller than $\ell(L)^{-C_4\beta}$ and they are negligible if $C_4$ is chosen large enough. It remains to upper bound the difference $\pi^{(1)}(\eta_x) - \pi^{(2)}(\eta_x)$. One has by translation invariance

$$\pi^{(1)}(\eta_x) = \pi^\lambda_{a+\ell(L)+1}(\eta_{x+a} = 1 | \eta_y = 0, \, 1 \leq y \leq a).$$

Recalling Theorem 1, this can be rewritten as

$$\pi^{(1)}(\eta_x) = \frac{\sum_{0 < y \le x} K_\lambda(a+y)\mathbf{P}_\lambda(x-y)\mathbf{P}_\lambda(\ell-x+1)}{\sum_{0 < y \le \ell+1} K_\lambda(a+y)\mathbf{P}_\lambda(\ell-y+1)}$$

where with some abuse of notation we let $\mathbf{P}_\lambda(n) := \mathbf{P}_\lambda(n \in \hat{S})$ and we wrote $\ell = \ell(L)$ for lightness of notation. Analogously, one can write

$$\pi^{(2)}(\eta_x) = \frac{\sum_{0 < y \le x} K_\lambda(a+y)\mathbf{P}_\lambda(x-y)\mathbf{R}_\lambda(\ell-x)}{\mathbf{P}_\lambda(\hat{S}_1 = \infty) + \sum_{0 < y \le \ell} K_\lambda(a+y)\mathbf{R}_\lambda(\ell-y)},$$

where $\mathbf{R}_\lambda(n) := \mathbf{P}_\lambda(\nexists i \ge 0, j > n : \hat{S}_i = j)$ (i.e. the probability that the transient renewal $\hat{S}$ contains no points after $n$) and altogether

$$\pi^{(1)}(\eta_x) - \pi^{(2)}(\eta_x) = \sum_{0 < y \le x} K_\lambda(a+y)\mathbf{P}_\lambda(x-y)$$

$$\times \left[ \frac{\mathbf{P}_\lambda(\ell-x+1)}{\sum_{0 < y \le \ell+1} K_\lambda(a+y)\mathbf{P}_\lambda(\ell-y+1)} - \frac{\mathbf{R}_\lambda(\ell-x)}{1-\lambda + \sum_{0 < y \le \ell} K_\lambda(a+y)\mathbf{R}_\lambda(\ell-y)} \right]$$
(6.15)

where we used $\mathbf{P}_\lambda(\hat{S}_1 = \infty) = 1 - \lambda$. Since the renewal $\hat{S}$ is transient, it is well known [10, App. A] that $\mathbf{P}_\lambda(n) \sim \frac{1}{(1-\lambda)^2}K_\lambda(n)$ for $n \to \infty$. As a consequence, one sees immediately that

$$|1 - \mathbf{R}_\lambda(n)| = O(n^{-\rho}).$$
(6.16)

It follows that

$$\frac{\mathbf{R}_\lambda(\ell-x)}{1-\lambda + \sum_{0 < y \le \ell} K_\lambda(a+y)\mathbf{R}_\lambda(\ell-y)} = \frac{1}{1-\lambda \sum_{y=1}^{a} K(y)} + O(\ell^{-\rho})$$

because the denominator of the term on the left hand side is equal to

$$1 - \sum_{y=1}^{a} K_\lambda(y) - \sum_{y=a+1}^{a+\ell} K_\lambda(y)[1 - \mathbf{R}_\lambda(\ell-(y-a))] - \sum_{y=a+\ell+1}^{\infty} K_\lambda(y)$$

and, decomposing the second sum according to $y \le a + \ell/2$ and $y > a + \ell/2$, we see this equals $1 - \sum_{y=1}^{a} K_\lambda(y) + O(\ell^{-\rho})$.

Now we need the following more accurate estimate, which is proven in Appendix B: there exists $w = w(\rho) \in (0, \rho)$ such that

$$\mathbf{P}_\lambda(n) = \frac{\lambda}{(1-\lambda)^2}K(n)(1 + O(n^{-w})).$$
(6.17)

To upper bound the first term in the square bracket in (6.15), we estimate separately numerator and denominator. The numerator is

$$\frac{\lambda}{(1-\lambda)^2}K(\ell-x+1)(1 + O(\ell^{-w}))$$

where we used the fact that $(\ell - x) \ge \ell/2$.

As for the denominator, it equals

$$\lambda \sum_{y=1}^{\ell/2} K(a+y)\mathbf{P}_\lambda(\ell-y+1) + \lambda \sum_{y=\ell/2+1}^{\ell+1} K(a+y)\mathbf{P}_\lambda(\ell-y+1).$$
(6.18)

By using (6.17) and the decreasing character of $K(\cdot)$, the first sum can be bounded below by

$$\frac{\lambda^2}{(1-\lambda)^2} K(\ell) \left\{ \sum_{y=a+1}^{a+\ell/2} K(y) \right\} (1 + O(\ell^{-w}))$$

$$= \frac{\lambda^2}{(1-\lambda)^2} K(\ell) \left\{ \sum_{y>a+1} K(y) + O(\ell^{-\rho}) \right\} (1 + O(\ell^{-w}))$$

$$= \frac{\lambda^2}{(1-\lambda)^2} K(\ell) \left\{ \sum_{y>a+1} K(y) \right\} (1 + O(\ell^{-w}))$$

and the second sum by $\lambda K(\ell + a + 1) \sum_{y=0}^{\ell/2} \mathbf{P}_\lambda(y)$.

Note that $\sum_{y=0}^{\infty} \mathbf{P}_\lambda(y)$ is the average total number of renewal points in $\hat{S}$: a simple computation on geometric random variables shows that it equals $1/(1-\lambda)$. From this fact, it is easy to deduce that $\sum_{y=0}^{\ell/2} \mathbf{P}_\lambda(y) = 1/(1-\lambda) + O(\ell^{-\rho})$. Observe also that, since we are assuming $a < L^g$, by taking $g$ sufficiently close to $\beta/\rho$ and $\beta$ sufficiently small one has $K(a + \ell + 1) = K(\ell)(1 + O(\ell^{-\rho}))$. Putting everything together, the quantity in (6.18) is lower bounded by

$$\lambda/(1-\lambda)^2 K(\ell) \left[ 1 - \lambda \sum_{y=1}^{a} K(y) \right] (1 + O(\ell^{-w})).$$

In conclusion, the l.h.s. of (6.15) is given by

$$\pi^{(1)}(\eta_x) - \pi^{(2)}(\eta_x)$$

$$= \frac{\sum_{0 < y \le x} K_\lambda(a+y) \mathbf{P}_\lambda(x-y)}{1 - \lambda \sum_{y=1}^{a} K(y)} \left[ \frac{K(\ell - x + 1)}{K(\ell)} - 1 + O(\ell^{-w}) \right]$$

$$\le (1-\lambda)^{-1} \sum_{0 < y \le x} K_\lambda(y) \mathbf{P}_\lambda(x-y) \left[ \frac{K(\ell - x + 1)}{K(\ell)} - 1 + O(\ell^{-w}) \right]$$

where we used the fact that $K(\cdot)$ is decreasing and $\sum_n K(n) = 1$.

Uniformly in the range $1 \le x \le \ell/2$ one has, for a positive constant $C(\rho)$,

$$\frac{K(\ell - x + 1)}{K(\ell)} - 1 \le C(\rho) \frac{x}{\ell}.$$

Also, the renewal equation gives

$$\sum_{0 < y \le x} K_\lambda(y) \mathbf{P}_\lambda(x-y) = \mathbf{P}_\lambda(x) \tag{6.19}$$

so that finally

$$\pi^{(1)}(\eta_x) - \pi^{(2)}(\eta_x) \le \frac{c(\lambda, \rho)}{x^{1+\rho}} \left( \frac{x}{\ell} + \frac{1}{\ell^w} \right). \tag{6.20}$$

Summing on $x \le \max(L^g - a, 1) \le \ell$, one gets a quantity which is $O(\ell^{-w}) = o(L^{-\beta})$ if $\beta$ is small enough. In particular, for $L$ large enough

$$\sum_{x=1}^{\max(L^g-a,1)} \left( \mu_t^{\sigma}(\eta_x) - \mu_t^{\sigma'}(\eta_x) \right) < \frac{1}{8L^\beta}. \tag{6.21}$$

Putting together Eqs. (6.13), (6.14) and (6.21) one obtains that the sum in (6.11) is smaller than $1/L^{\beta}$, which finishes the proof of Proposition 3.

$\square$

**Remark 4.** *It is easy to see that, for $\lambda$ sufficiently small, the mixing time is actually $O(\log L)$. For this, it is sufficient to check that if $\lambda \ll 1$ the dynamics contracts the Hamming distance, i.e. (6.2) holds for some $t = O(1)$ and $\gamma > 0$, whenever $d(\sigma, \sigma') = 1$. Then, by Theorem 6, Equation (6.10) and the fact that the Hamming distance is upper bounded by the system size, one gets that $\|\mu_t^{\sigma} - \mu_t^{\sigma'}\| < 1/(2e)$ for every couple initial conditions $(\sigma, \sigma')$, for some $t = O(\log |\Lambda|)$, which implies the claim.*

6.2. **Relaxation time.** We prove here the second part of Theorem 3, i.e. $T_{rel} \leq C_4(\lambda)$. The proof is decomposed in two steps. First, in Lemma 9 we show that the Dirichlet form of the initial dynamics is comparable to the Dirichlet form of a block dynamics whenever the size of the blocks remains finite, i.e depending on $\lambda$ but not on the size $L$ of the big box. Thus, it is sufficient to prove that the relaxation time of the block dynamics does not depend on $L$. In the second step (Lemma 10) we show, by using the path coupling argument by Bubley and Dyer [4], that if the size of these blocks is sufficiently large then the relaxation time of the block dynamics is independent of $L$.

The first step would be standard (cf. e.g. [13, Sec. 3]) for, say, the heat bath dynamics of the nearest neighbor Ising model, where transition rates $c_x(\sigma, \sigma')$ are uniformly bounded from below and depend on $\sigma, \sigma'$ only in a finite neighborhood of the site $x$ where the Poisson clock rings (finite-range interactions). In our case, both uniform positivity of the rates and finite-range dependence fails, cf. the definitions (2.10)-(2.11) of creation/destruction rates.

The block dynamics is denoted by $\{\hat{\eta}(t)\}_{t \geq 0}$. The size of the blocks equals $2\ell + 1$, where in this section $\ell$ depends on $\lambda$ but not on $L$. For any $x \in \{1, \ldots, L-1\}$ let $\Lambda_{\ell}(x) = \{y : |y - x| \leq \ell\} \cap \{1, \ldots, L-1\}$ and $\Lambda_{\ell}^c(x)$ the complementary set of $\Lambda_{\ell}(x)$ in $\{0, \ldots, L\}$. The generator $\hat{\mathcal{L}}$ of $\{\hat{\eta}(t)\}_{t \geq 0}$ is given by

$$(\hat{\mathcal{L}}f)(\eta) = \sum_{x=1}^{L-1} \left[ \hat{Q}_x f(\eta) - f(\eta) \right] = \sum_{x=1}^{L-1} \sum_{\tilde{\eta}, \eta} \hat{c}_x(\eta, \tilde{\eta}) \left[ f(\tilde{\eta}) - f(\eta) \right]$$

where $(\hat{Q}_x f)(\eta) = \pi_L^{\lambda} \left( f \,|\, \eta_{\Lambda_{\ell}^c(x)} \right)$. The rate $\hat{c}_x(\eta, \tilde{\eta})$ is zero if the restrictions to $\Lambda_{\ell}^c(x)$ of $\eta$ and $\tilde{\eta}$ are distinct. If instead $\eta_{\Lambda_{\ell}^c(x)} = \tilde{\eta}_{\Lambda_{\ell}^c(x)}$, let $a := a(x, \eta)$ (resp. $b := b(x, \eta)$) be the rightmost particle of $\eta$ in $\Lambda_{\ell}^c(x)$ to the left of $x$ (resp. the leftmost particle of $\eta$ in $\Lambda_{\ell}^c(x)$ to the right of $x$). Then,

$$\hat{c}_x(\eta, \tilde{\eta}) = \frac{W_x^{a,b}(\tilde{\eta}|_{\Lambda_{\ell}(x)})}{\sum_{\sigma \in \{0,1\}^{\Lambda_{\ell}(x)}} W_x^{a,b}(\sigma)} \tag{6.22}$$

with $W_x^{a,b}(\sigma) = \lambda^k \prod_{j=0}^{k} K(z_{j+1} - z_j)$, $z_0 = a$, $z_{k+1} = b$, and $z_j, j = 1, \ldots, k$ the ordered points $z$ of $\Lambda_{\ell}(x)$ for which $\sigma_z = 1$. In this formula, it is understood that for $k = 0$, $W_x^{a,b}(\sigma) = K(b - a)$. If $\omega \in \Omega_L$, we use the short notation $W_x^{a,b}(\omega)$ for $W_x^{a,b}(\omega_{\Lambda_{\ell}(x)})$.

The Dirichlet form associated to $\hat{\mathcal{L}}$ is given by

$$\hat{\mathcal{E}}_L(f;f) = \pi_L^\lambda \left[ f \, ; \, -\hat{\mathcal{L}}f \right] = \frac{1}{2} \sum_{x=1}^{L-1} \sum_{\eta,\tilde{\eta}} \pi_L^\lambda(\eta) \, \hat{c}_x(\eta,\tilde{\eta}) \, [f(\tilde{\eta}) - f(\eta)]^2 \, .$$

**Lemma 9.** *For any $\ell \geq 1$, there exists a constant $C := C(\ell, \lambda)$ independent of $L$ such that for any function $f : \Omega_L \to \mathbb{R}$,*

$$\hat{\mathcal{E}}_L(f;f) \leq C \, \mathcal{E}_L(f;f).$$

*Proof.* In the proof, $C$ denotes a positive constant (depending only on $\lambda, \rho$ and $\ell$, but not on $L$) which is not the same at each occurrence.

For any configurations $\eta, \tilde{\eta}$ which coincide outside of $\Lambda_\ell(x)$, we denote by $\gamma_x(\eta, \tilde{\eta}) = \left\{ \eta = \eta^{(0)}, \eta^{(1)}, \dots, \eta^{(p)} = \tilde{\eta} \right\}$ a sequence of configurations such that $\eta^{(i+1)}$ is obtained from $\eta^{(i)}$ by adding or deleting a particle at some site $x_i \in \Lambda_\ell(x)$. The sequence is chosen such that $p := p(\eta, \tilde{\eta})$ is minimal. Observe that $p$ is at most $2\ell + 1$, independently of $\eta, \tilde{\eta}, x, L$.

Write $f(\tilde{\eta}) - f(\eta) = \sum_{j=0}^{p-1} (f(\eta^{(j+1)}) - f(\eta^{(j)}))$ and use Schwarz inequality:

$$\hat{\mathcal{E}}_L(f;f) \leq C \sum_x \sum_{\eta,\tilde{\eta}} \pi_L^\lambda(\eta) \, \hat{c}_x(\eta,\tilde{\eta}) \sum_{j=0}^{p-1} \left[ f(\eta^{(j+1)}) - f(\eta^{(j)}) \right]^2$$

$$= C \sum_{x,\eta,\tilde{\eta}} \sum_{j=0}^{p-1} \frac{\pi_L^\lambda(\eta)}{\pi_L^\lambda(\eta^{(j)})} \frac{\hat{c}_x(\eta,\tilde{\eta})}{c(\eta^{(j)},\eta^{(j+1)})} c(\eta^{(j)},\eta^{(j+1)}) \left[ f(\eta^{(j+1)}) - f(\eta^{(j)}) \right]^2 \pi_L^\lambda(\eta^{(j)}).$$

Here, $c(\eta^{(j)}, \eta^{(j+1)})$ stands for the rate of change from $\eta^{(j)}$ to $\eta^{(j+1)}$ for the original dynamics $\{\eta_t\}_t$. By inverting the different sums, we are left to show that for any $z \in \{1, \dots, L-1\}$, any $\omega \in \Omega_L$,

$$\sum_{x \in \Lambda_\ell(z)} \sum_{(\eta,\tilde{\eta})} R(\omega, z, x, \eta, \tilde{\eta}) = \sum_{x \in \Lambda_\ell(z)} \sum_{(\eta,\tilde{\eta})} \mathbf{1}_{\{(\omega,\omega^z) \in \gamma_x(\eta,\tilde{\eta})\}} \, \frac{\pi_L^\lambda(\eta)}{\pi_L^\lambda(\omega)} \frac{\hat{c}_x(\eta,\tilde{\eta})}{c(\omega,\omega^z)}$$

is bounded by a constant independent of $\omega$, $z$ and $L$. Observe that the two sums, because of the indicator function, can be carried over a set whose cardinal is bounded by a constant (independent of $\omega, z, L$). Therefore, it is sufficient to bound each term of the sum. Since $(\omega, \omega^z) \in \gamma_x(\eta, \tilde{\eta})$, the restrictions to $\Lambda_\ell^c(x)$ of $\omega, \omega^z, \eta, \tilde{\eta}$ are equal and $a(x, \omega) = a(x, \eta) = a(x, \tilde{\eta})$ and similarly for $b$. Thus,

$$R(\omega, z, x, \eta, \tilde{\eta}) = \mathbf{1}_{\{(\omega,\omega^z) \in \gamma_x(\eta,\tilde{\eta})\}} \frac{1}{c(\omega, \omega^z)} \frac{W_x^{a,b}(\eta) W_x^{a,b}(\tilde{\eta})}{W_x^{a,b}(\omega) \sum_{\sigma \in \{0,1\}^{\Lambda_\ell(x)}} W_x^{a,b}(\sigma)}. \quad (6.23)$$

To fix ideas, we assume that $\omega_z = 0$ (the other case can be treated similarly). The rate $c(\omega, \omega^z)$ corresponds then to the creation of a particle at site $z$. Since $(\omega, \omega^z) \in \gamma_x(\eta, \tilde{\eta})$ we have $\eta_z = 0$, $\tilde{\eta}_z = 1$. Let $u \geq a$ (resp. $v \leq b$) be the rightmost particle of $\omega$ to the left of $z$ (resp. the leftmost one to the right of $z$). We have from (2.10)

$$c(\omega, \omega^z) = \frac{\lambda K(v-z) K(z-u)}{\lambda K(v-z) K(z-u) + K(v-u)}. \quad (6.24)$$

We distinguish two cases according to the positions of $u, v$ with respect to $a, b$.

First case: $a < u$ (the case $v < b$ being very similar). This means that $u \in \Lambda_\ell(x)$, so that $K(z-u) \geq K(2\ell+1)$. Then, from (6.24) one realizes immediately that $c(\omega, \omega^z)$ is lower bounded by some positive constant. On the other hand, one has that either $\eta_u = 1$ or $\tilde{\eta}_u = 1$ (or both) since the path $\gamma_x(\eta, \tilde{\eta})$ is assumed to be of minimal length. Assume for definiteness that the first case occurs. Then, it is easy to see that $W_x^{a,b}(\eta)/W_x^{a,b}(\omega) \leq C$ and then (6.23) is upper bounded by a positive constant (clearly $\sum_{\sigma \in \{0,1\}^{\Lambda_\ell(x)}} W_x^{a,b}(\sigma) \geq W_x^{a,b}(\tilde{\eta})$).

Second case: $u = a$ and $b = v$. This means that the restriction of $\omega$ to $\Lambda_\ell(x)$ is empty. Since $(\omega, \omega^z) \in \gamma_x(\eta, \tilde{\eta})$, $\tilde{\eta}_z = 1$, and we get that $W_x^{a,b}(\tilde{\eta}) \leq C K(z - a) K(b - z)$. It follows easily (using the definition of $K(\cdot)$) that

$$\frac{W_x^{a,b}(\tilde{\eta})}{c(\omega, \omega^z) W_x^{a,b}(\omega)} \leq C \left[ \frac{K(z-a)K(b-z)}{K(b-a)} + \frac{1}{\lambda} \right] \leq C'$$

so that $R(\omega, z, x, \eta, \tilde{\eta}) \leq C$ also in this case. $\qquad \square$

**Lemma 10.** *There exists $\ell := \ell(\lambda)$ such that the block dynamics $\{\hat{\eta}(t)\}_{t \geq 0}$ has a spectral gap bounded below by a constant $C := C(\ell, \lambda) > 0$ independent of $L$.*

*Proof.* We recall that the Hamming distance $d(\sigma, \omega)$ between two configurations $\sigma, \omega \in \Omega_L$ is defined by $d(\sigma, \omega) = \sum_{x=1}^{L-1} |\sigma_x - \omega_x|$. By [4] (see also Theorem 6 above), it is sufficient to find a generator $\bar{\mathcal{L}}$ acting on functions $f : \Omega_L \times \Omega_L \ni (\eta, \eta') \mapsto \mathbb{R}$ which acts like $\mathcal{L}$ if $f$ depends only on $\eta$ (or only on $\eta'$) and such that there exists a positive constant $\delta > 0$ satisfying

$$(\bar{\mathcal{L}}d)(\omega, \omega^z) \leq -\delta \qquad (6.25)$$

for every $\omega \in \Omega_L$ and $z \in \{1, \ldots, L-1\}$. We write

$$(\bar{\mathcal{L}}f)(\sigma, \omega) = \sum_{x=1}^{L-1} [\nu_x^{\sigma, \omega}(f) - f(\sigma, \omega)]$$

where $\nu_x^{\sigma, \omega}$ is some coupling between $\pi_L^\lambda(\cdot|\sigma|_{\Lambda_\ell^c(x)})$ and $\pi_L^\lambda(\cdot|\omega|_{\Lambda_\ell^c(x)})$. Let $\sigma \in \Omega_L$ and $z \in \{1, \ldots, L-1\}$ such that $\sigma_z = 0$. We know that $\pi_L^\lambda(\cdot|\sigma^z|_{\Lambda_\ell^c(x)})$ dominates stochastically $\pi_L^\lambda(\cdot|\sigma|_{\Lambda_\ell^c(x)})$ and in this case we require $\nu_x^{\sigma, \sigma^z}$ to be a monotone coupling. We denote by $0 \leq a < z$ (resp. $z < b \leq L$) the rightmost site to the left of $z$ (resp. the leftmost site to the right of $z$) which is occupied in the configuration $\sigma$. We have (using $d(\sigma, \sigma^z) = 1$)

$$(\bar{\mathcal{L}}d)(\sigma, \sigma^z) = \sum_{x:|x-z| \leq \ell} \left[ \nu_x^{\sigma, \sigma^z}(d) - 1 \right] + \sum_{x:|x-z| > \ell} \left[ \nu_x^{\sigma, \sigma^z}(d) - 1 \right].$$

If $|x - z| \leq \ell$ then $\sigma|_{\Lambda_\ell^c(x)} = \sigma^z|_{\Lambda_\ell^c(x)}$ and we have $\nu_x^{\sigma, \sigma^z}(d) = 0$. Thus the first sum is equal to $-|\{x \in \{1, \ldots, L-1\}; |x-z| \leq \ell\}| \leq -\ell$. In the second sum, observe that if $x > b + \ell$ (or similarly if $x < a - \ell$) then after the update of the box $\Lambda_\ell(x)$, the distance between $\sigma$ and $\sigma^z$ is again 1, so that these terms do not contribute. The reason is that the presence of a particle in $b$ "screens" the effect of the discrepancy in $z$: one has $\pi_L^\lambda(\cdot|\sigma|_{\Lambda_\ell^c(x)})$ and $\pi_L^\lambda(\cdot|\sigma^z|_{\Lambda_\ell^c(x)})$, so the coupling is diagonal. What is left of the second sum is

$$\sum_{x \in I_\ell(a,z,b)} \sum_{y \in \Lambda_\ell(x)} \left\{ \pi_L^\lambda \left[ \eta_y \,\Big|\, \eta|_{\Lambda_\ell^c(x)} = \sigma^z|_{\Lambda_\ell^c(x)} \right] - \pi_L^\lambda \left[ \eta_y \,\Big|\, \eta|_{\Lambda_\ell^c(x)} = \sigma|_{\Lambda_\ell^c(x)} \right] \right\}$$

$$\qquad (6.26)$$

where $I_\ell(a, z, b) = \{x \in \{1, \ldots, L-1\} : |x - z| > \ell, a - \ell \leq x \leq b + \ell\}$. It remains to show there exists a constant $\ell_0 := \ell_0(\lambda, \rho)$ such that if $\ell \geq \ell_0$ then this sum is less than $\ell/4$. We restrict to the sum carried over the sites $z + \ell < x \leq b + \ell$ since the other case can be treated similarly.

By monotonicity, we can and will assume that $\sigma$ is empty in $\{1, \ldots, z-1\}$. This way, the first term in (6.26) does not change and the second one decreases, so altogether the sum increases.

For any $\alpha < \beta$, we denote by $\mu_x^{\alpha, \beta}$ the law of the position $s$ of the leftmost particle in $\Lambda_\ell(x) \subset \{\alpha, \ldots, \beta\}$ conditionally to the fact that there is a particle on site $\alpha$ and on site $\beta$ and no particles in $\{\alpha + 1, \ldots \beta - 1\} \cap \Lambda_\ell^c(x)$. Note that there could be zero particles in $\Lambda_\ell(x) \subset \{\alpha, \ldots, \beta\}$, i.e., $\mu_x^{\alpha, \beta}$ is a sub-probability.

For any $x$ such that $z + \ell < x \leq b + \ell$, we denote $b_x$ the leftmost site $y > (x + \ell) \wedge L$ such that $\sigma_y = 1$ (observe that $b_x = b$ if $x \leq b - \ell - 1$). Then we have

$$\sum_{y \in \Lambda_\ell(x)} \left\{ \pi_L^\lambda \left[ \eta_y \,\Big|\, \eta|_{\Lambda_\ell^c(x)} = \sigma^z|_{\Lambda_\ell^c(x)} \right] - \pi_L^\lambda \left[ \eta_y \,\Big|\, \eta|_{\Lambda_\ell^c(x)} = \sigma|_{\Lambda_\ell^c(x)} \right] \right\}$$

$$= \sum_{s \in \Lambda_\ell(x)} \left( \mu_x^{z, b_x}(s) - \mu_x^{0, b_x}(s) \right) \pi_{[s, b_x]}^\lambda \left( \sum_{y=s}^{x+\ell} \eta_y \,\Big|\, \eta_u = 0 \text{ for } x + \ell < u < b_x \right) \quad (6.27)$$

$$\leq C \min \left\{ \|\mu_x^{z, b_x} - \mu_x^{0, b_x}\|, \sum_{s \in \Lambda_\ell(x)} \mu_x^{z, b_x}(s) \right\}$$

where the last inequality follows because, in the delocalized phase, the average number of particles in a box is uniformly upper bounded by a constant. $\pi_{[s, b_x]}^\lambda$ denotes the equilibrium measure $\pi_{b_x - s + 1}^\lambda$ translated by $s$. Below we estimate the last line in (6.27) distinguish various cases according to the values of $x, b_x$ etc. Recall that $\mu_x^{z, b_x}$ is a sub-probability, so in many cases (when the box $\Lambda_\ell(x)$ is sufficiently far from $z, b_x$) the sum $\sum_{s \in \Lambda_\ell(x)} \mu_x^{z, b_x}(s)$ is quite small. In the remaining cases, we really have to estimate the variation distance between $\mu_x^{z, b_x}$ and $\mu_x^{0, b_x}$.

By using the same notations as in the proof of Proposition 3, we can write

$$\mu_x^{\alpha, \beta}(s) = \frac{\sum_{s \leq v \leq x + \ell} K_\lambda(s - \alpha) \mathbf{P}_\lambda(v - s) K_\lambda(\beta - v)}{K_\lambda(\beta - \alpha) + \sum_{x - \ell \leq u \leq x + \ell} K_\lambda(u - \alpha) \mathbf{P}_\lambda(v - u) K_\lambda(\beta - v)}.$$

**Case 1:** $z + \ell + 1 \leq x \leq b - \ell - 1$. We have $b_x = b$. By using the renewal equation (6.19) and $\mathbf{P}_\lambda(n) \leq CK(n)$ (cf. (6.17)), we get

$$\sum_{x=z+\ell+1}^{b-\ell-1} \sum_{s \in \Lambda_\ell(x)} \mu_x^{z, b_x}(s) \leq \sum_{x=z+\ell+1}^{b-\ell-1} \sum_{s \in \Lambda_\ell(x)} \frac{\sum_{s \leq v \leq b} K_\lambda(s - z) \mathbf{P}_\lambda(v - s) K_\lambda(b - v)}{K_\lambda(b - z)}$$

$$\leq C \sum_{x=z+\ell+1}^{b-\ell-1} \sum_{s \in \Lambda_\ell(x)} \frac{K(s - z) K(b - s)}{K(b - z)}.$$

The latter sum can be written as

$$\sum_{s=z+1}^{b-1} \left[ (s - z)^{-1} + (b - s)^{-1} \right]^{1+\rho} N_\ell(z, s, b)$$

with $N_\ell(z, s, b)$ the number of sites $x$ such that $s \in \Lambda_\ell(x)$:

$$N_\ell(z, s, b) = \begin{cases} s - z \text{ if } z + 1 \le s \le z + 2\ell, \\ 2\ell + 1 \text{ if } z + 2\ell + 1 \le s \le b - 2\ell - 1, \\ b - s \text{ if } b - 2\ell \le s \le b - 1. \end{cases}$$

By using $\left[ (s - z)^{-1} + (b - s)^{-1} \right]^{1+\rho} \le C \left[ (s - z)^{-(1+\rho)} + (b - s)^{-(1+\rho)} \right]$, we are left to estimate

$$\sum_{s=z+1}^{z+2\ell} \frac{1}{(s-z)^\rho} + (2\ell + 1) \sum_{s=z+2\ell+1}^{b-2\ell-1} \frac{1}{(s-z)^{1+\rho}} + \sum_{s=b-2\ell}^{b-1} \frac{b-s}{(s-z)^{1+\rho}}$$

and

$$\sum_{s=z+1}^{z+2\ell} \frac{(s-z)}{(b-s)^{1+\rho}} + (2\ell + 1) \sum_{s=z+2\ell+1}^{b-2\ell-1} \frac{1}{(b-s)^{1+\rho}} + \sum_{s=b-2\ell}^{b-1} \frac{1}{(b-s)^\rho}$$

which are of order $\ell^{1-\rho}$ uniformly in $z$ and $b$.

**Case 2:** $\sup(b - \ell, z + \ell + 1) \le x \le b + \ell$ (note that there are at most $2\ell + 1$ possible values for $x$). Let $\varepsilon > 0$ sufficiently small and assume to simplify notations that $\varepsilon\ell \ge 2$ is an integer. In the sequel, the constant $C_\varepsilon$ depends on $\varepsilon$ while $C$ does not.

**Case 2.a:** $z < x - \ell \le z + \varepsilon\ell$: As in Case 1, we have

$$\sum_{s \in \Lambda_\ell(x)} \mu_x^{z, b_x}(s) \le C \sum_{s \in \Lambda_\ell(x)} \frac{K(s-z)K(b_x - s)}{K(b_x - z)} = C \sum_{s \in \Lambda_\ell(x)} \left[ (b_x - s)^{-1} + (s - z)^{-1} \right]^{1+\rho}$$

which is finite because $\sum_{n \ge 1} n^{-(1+\rho)} < +\infty$. Since there are at most $\varepsilon\ell$ possible values of $x$, the contribution to (6.26) is upper bounded by $C\varepsilon\ell$.

**Case 2.b:** $x - \ell > z + \varepsilon\ell$ and $b_x - (x + \ell) \ge \varepsilon\ell$. Again

$$\sum_{s \in \Lambda_\ell(x)} \mu_x^{z, b_x}(s) \le C \sum_{s \in \Lambda_\ell(x)} \left[ (b_x - s)^{-1} + (s - z)^{-1} \right]^{1+\rho}$$

$$\le C \sum_{s \in \Lambda_\ell(x)} \left[ (x + (1 + \varepsilon)\ell - s)^{-(1+\rho)} + (s - z)^{-(1+\rho)} \right].$$

It is easy to show that summing over the allowed $x$ this term contribute to (6.26) at most $C_\varepsilon \ell^{(1-\rho)}$.

**Case 2.c:** $x - \ell > z + \varepsilon\ell$ and $b_x - (x + \ell) < \varepsilon\ell$: Here at last we use the upper bound $\|\mu_x^{z, b_x} - \mu_x^{0, b_x}\|$ in (6.27). We define $F(\alpha, s) = N(\alpha, s)/D(\alpha, s)$ where

$$N(\alpha, s) = \left( 1 - \frac{b_x - s}{b_x - \alpha} \right)^{-(1+\rho)} \sum_{s \le v \le x + \ell} \mathbf{P}_\lambda(v - s) K_\lambda(b_x - v)$$

$$D(\alpha, s) = 1 + \sum_{x - \ell \le u \le v \le x + \ell} \mathbf{P}_\lambda(v - u) K_\lambda(b_x - v) \left( 1 - \frac{b_x - u}{b_x - \alpha} \right)^{-(1+\rho)}.$$

Note that $\mu_x^{z, b_x}(s) = F(z, s)$, $\mu_x^{0, b_x}(s) = F(0, s)$ so that

$$\left| \mu_x^{z, b_x}(s) - \mu_x^{0, b_x}(s) \right| \le \int_0^z |\partial_\alpha F(\alpha, s)| \, d\alpha.$$

Observe now (using also $D \geq 1$) that $|\partial_\alpha F| \leq |\partial_\alpha N| + N|\partial_\alpha D|$. Using $\mathbf{P}_\lambda(n) \leq CK(n)$ and the renewal equation (6.19) we get easily

$$N \leq C \left( 1 - \frac{b_x - s}{b_x - \alpha} \right)^{-(1+\rho)} \frac{1}{(b_x - s)^{(1+\rho)}} \leq C_\varepsilon \frac{1}{(b_x - s)^{(1+\rho)}}$$

$$|\partial_\alpha N| \leq C \left( 1 - \frac{b_x - s}{b_x - \alpha} \right)^{-(2+\rho)} \frac{1}{(b_x - \alpha)^2(b_x - s)^\rho} \leq C_\varepsilon \frac{1}{(b_x - \alpha)^2(b_x - s)^\rho}.$$

The last inequalities follows from $\alpha \leq z$, $b_x - s = b_x - s + (s - z) \geq b_x - s + \varepsilon\ell$, the fact that the function $x \mapsto x/(x + \varepsilon\ell)$ is increasing and $b_x - s \leq (2 + \varepsilon)\ell$. Using similar bounds one obtains

$$|\partial_\alpha D| \leq C_\varepsilon \frac{1}{(b_x - \alpha)^2} \sum_{u \in \Lambda_\ell(x)} \frac{1}{(b_x - u)^\rho} \leq C_\varepsilon \frac{\ell^{1-\rho}}{(b_x - \alpha)^2}.$$

Thus, integrating w.r.t. $\alpha$ and using $b_x - z \geq 2\varepsilon\ell$, we obtain

$$\int_0^z |\partial_\alpha F(\alpha, s)| \, d\alpha \leq \frac{C(\varepsilon)}{\ell} \left[ \frac{1}{(b_x - s)^\rho} + \frac{\ell^{1-\rho}}{(b_x - s)^{1+\rho}} \right].$$

The sum over $s \in \Lambda_\ell(x)$ of the r.h.s. is bounded above by $C_\varepsilon \ell^{-\rho}$. Since we have now to take the sum over a set of $x$'s of cardinality at most $2\ell + 1$, we conclude that the contribution of these $x$'s in (6.26) is bounded by $C_\varepsilon \ell^{(1-\rho)}$.

Putting everything together, we have proven that (6.26) is upper bounded by $C_\varepsilon \ell^{1-\rho} + C\varepsilon\ell$. Taking $\varepsilon$ sufficiently small and then $\ell$ sufficiently large, this is smaller than $\ell/4$, which concludes the proof. □

## 7. The critical point

Here $\lambda = 1$; we write $\pi_L := \pi_L^{\lambda=1}$ and (see Theorem 1)

$$Z_L := Z_L(\lambda = 1) \sim C_1 \times \frac{1}{L^{1-\rho}}, \quad C_1 > 0.$$

It would be easy to prove $\mathrm{T_{mix}} \geq C_5 L^\rho$ as we did in the localized phase: starting from the empty configuration, for times much shorter than $L^\rho$ the chance of creating a particle between $\varepsilon L$ and $(1 - \varepsilon)L$ is small. On the other hand, at equilibrium

$$\pi_L(\eta_x = 0 \text{ for every } \varepsilon L \leq x \leq (1 - \varepsilon)L) \sim C_2 \times \varepsilon^{2\rho}, \quad C_2 > 0$$

and the mixing time lower bound follows for $\varepsilon$ sufficiently small.

What we do instead is to prove that $C_5 L^\rho \leq \mathrm{T_{rel}} \leq C_6 L^{1+\rho}$, which implies directly the mixing time bounds via (2.16) and the fact that, as proven in Section 5.6 (cf. (5.23)), one has $\log(1/\pi^*) = O(L)$.

To get a lower bound for $\mathrm{T_{rel}}$, recall its variational definition (2.13) and choose $f(\eta) = \sum_{x=1}^{L-1} \eta_x$ (the number of particles between 1 and $L - 1$ in the configuration $\eta$). It is well known (and can be seen by an elementary computation) that $\pi_L(f) \sim c_1 L^\rho$ and $\pi_L(f; f) \sim c_2 L^{2\rho}$ for some $c_1, c_2 > 0$. Next we compute the Dirichlet form. By reversibility $(1 - \eta_x)c_x(\eta)\pi(\eta) = \pi(\eta^x)\eta_x^x d_x(\eta^x)$ so that for every $f$

$$\mathcal{E}_L(f; f) = \sum_{x=1}^{L-1} \sum_{\eta \in \Omega_L} \pi_L(\eta) \mathbf{1}_{\{\eta_x=1\}} d_x(\eta)(f(\eta) - f(\eta^x))^2.$$

In our case $(f(\eta) - f(\eta^x))^2 = 1$ and the destruction rates are upper bounded by a constant, so that

$$\mathcal{E}_L(f; f) \leq c\,\pi_L\left(\sum_{x=1}^{L-1} \mathbf{1}_{\eta_x = 1}\right) = c\,\pi_L(f) \leq c' L^\rho$$

with $c, c'$ some positive constants. Altogether, we have proven that the spectral gap of the chain is $O(L^{-\rho})$.

Finally, the upper bound $\mathrm{T}_{\mathrm{rel}} \leq C_6 L^{1+\rho}$. Note that the proof of Lemma 2 (which holds for every $\lambda$) gives directly $\mathrm{T}_{\mathrm{rel}} \leq c L^{2+\rho}$: we want to improve on this estimate, using the fact that we are at the critical point.

From (5.25), (5.26), the fact that $|\gamma(\eta, \eta')| \leq L$ one sees that to prove $\mathrm{T}_{\mathrm{rel}} = O(L^{1+\rho})$ it suffices to show there exists a positive constant $C$ such that, for every $\sigma \in \Omega_L$ and $z = 1, \ldots, L-1$ satisfying $\sigma_z = 1$,

$$Z_L^{-1} \sum_{\substack{\eta \in \Omega_L : \eta_z = \sigma_x, x \geq z \\ \eta' \in \Omega_L : \eta'_x = \sigma_x, x < z, \eta'_z = 0}} \frac{\prod_{j=1}^{b-1} K(x_{j+1}^\eta - x_j^\eta) \prod_{j=c}^{n(\eta')} K(x_{j+1}^{\eta'} - x_j^{\eta'})}{K(z - x_{a-1}^\sigma)} \leq C\,L^\rho \quad (7.1)$$

where we recall that we identify the configuration $\sigma$ with the positions of its particles up to position $L$: $0 = x_0^\sigma < \cdots < x_{n(\sigma)+1}^\sigma = L$ and similarly for $\eta, \eta'$ (we also let $a$ the index such that $x_a^\sigma = z$, $b$ the index such that $x_b^\eta = z$ and $c$ the largest index such that $x_c^{\eta'} < z$). One has, since $x_b^\eta = z$,

$$\sum_{\eta \in \Omega_L : \eta_x = \sigma_x, x \geq z} \prod_{j=1}^{b-1} K(x_{j+1}^\eta - x_j^\eta) = Z_z \sim c\,z^{\rho-1}. \quad (7.2)$$

Next, since as we already observed that $x_c^{\eta'} = x_{a-1}^\sigma$, one has

$$\sum_{\substack{\eta' : \eta'_z = 0, \\ \eta'_x = \sigma_x, x < z}} \prod_{j=c}^{n(\eta')} K(x_{j+1}^{\eta'} - x_j^{\eta'}) = \mathbf{P}(L - x_{a-1}^\sigma \in S, S \cap \{1, \ldots, (z - x_{a-1}^\sigma)\} = \emptyset) \quad (7.3)$$

where like in (2.3) we let $\mathbf{P}$ denote the law of the renewal $S$ with inter-arrival distribution $K(\cdot)$. The probability in (7.3) equals

$$\sum_{y=(z - x_{a-1}^\sigma)+1}^{L} K(y) \mathbf{P}(L - y \in S).$$

Using the fact that $\mathbf{P}(L \in S) = Z_L(\lambda = 1) \sim c L^{\rho-1}$ one sees easily that

$$\sum_{y=(z - x_{a-1}^\sigma)+1}^{L} K(y) \mathbf{P}(L - y \in S) \leq c \times \begin{cases} \frac{L^{\rho-1}}{(z - x_{a-1}^\sigma)^\rho} & \text{if} \quad (z - x_{a-1}^\sigma) \leq L/2 \\ \frac{(L - z + x_{a-1}^\sigma)^\rho}{(z - x_{a-1}^\sigma)^{1+\rho}} & \text{if} \quad (z - x_{a-1}^\sigma) > L/2. \end{cases} \quad (7.4)$$

Putting together (7.2) and (7.4) one sees that the l.h.s. of (7.1) is upper bounded by a constant times $L^\rho$ as we wanted. $\qquad \square$

## Appendix A. Proof of Lemma 8

We first prove that $n^{-1} f(T(n, \ell))$ converges in law to a non-degenerate stable law of parameter $\rho$. The only point which requires some care is that the law $Q(\cdot)$ depends on the parameter $\ell$. This will prove the claim (5.16) for $n \geq n_0, \ell \geq \ell_0$ for some finite $n_0, \ell_0$. For every positive $\mu$ we have

$$\mathbb{E}\left(\exp\left\{-(\mu/n)\, f(T(n, \ell))\right\}\right) = \left[\mathbb{E}\left(\exp(-(\mu/n)\, f(1))\right)\right]^{T(n,\ell)} = \left[1 + (\mathfrak{L}(\mu/n) - 1)\right]^{T(n,\ell)} \tag{A.1}$$

where $\mathfrak{L}(\mu)$ is the Laplace transform of the law $Q$ defined in (5.15):

$$\mathfrak{L}(\mu) = \sum_{j=0}^{\infty} Q(j) e^{-\mu j} = 1 - \sum_{j=1}^{\infty} Q(j)(1 - e^{-\mu j}).$$

Let us start from the identity

$$Q(j) = c_0 \lambda K((j+1)\ell) - c_0 \frac{[\lambda K((j+1)\ell)]^2}{\lambda K((j+1)\ell) + 1} + \tilde{K}(j) \left[\prod_{i=j+1}^{\infty}(1 - \tilde{K}(i)) - 1\right]. \tag{A.2}$$

The dominant term is the first one: one has

$$\begin{aligned}
\sum_{j=1}^{\infty} K((j+1)\ell)(1 - e^{-\mu j}) &= \frac{\mathcal{C}_K \mu^\rho}{\ell^{1+\rho}} \mu \sum_{j=1}^{\infty} \frac{1 - e^{-\mu j}}{((j+1)\mu)^{1+\rho}} \\
&= \frac{\mathcal{C}_K \mu^\rho}{\ell^{1+\rho}} \left(\int_0^\infty \frac{1 - e^{-x}}{x^{1+\rho}} dx + o(\mu)\right)
\end{aligned} \tag{A.3}$$

as $\mu \to 0$. The second term in (A.2) gives a negligible contribution. Indeed, simply note that

$$\sum_{j=1}^{\infty} \left(K((j+1)\ell)\right)^2 \left(1 - e^{-\mu j}\right) = O(\mu \ell^{-2-2\rho}). \tag{A.4}$$

As for the last term, it is easy to see that for $\ell$ sufficiently large

$$\left|\prod_{i=j+1}^{\infty}(1 - \tilde{K}(i)) - 1\right| \leq \frac{c(\lambda)}{\ell^{1+\rho}(j+1)^\rho} \tag{A.5}$$

(just use that for $\ell$ large one has that $\tilde{K}(i)$ is small for every $i$, so that $(1 - \tilde{K}(i)) \approx \exp(-\tilde{K}(i)) \approx \exp(-c(\lambda)/(\ell^{1+\rho}(i+1)^{1+\rho})))$. As a consequence, one sees that the contribution to $\mathfrak{L}(\mu)$ from the last term in (A.2) is $O(\mu^{\min(1,2\rho)}\ell^{-2-2\rho})$. Altogether,

$$\mathfrak{L}(\mu) - 1 = -A \frac{\mu^\rho}{\ell^{1+\rho}} \left[1 + O\left(\frac{\mu^{\min((1-\rho),\rho)}}{\ell^{1+\rho}}\right) + o(\mu)\right] \tag{A.6}$$

for some positive constant $A$ depending only on $\lambda, \rho$ and $c_0$. Plugging this estimate into (A.1) and recalling that $T(n, \ell) = \lfloor n^\rho \ell^{1+\rho} \rfloor$ we obtain the desired convergence for $n, \ell \to \infty$.

Finally, we prove the claim (5.16) for $\ell \geq \ell_0$ and $n \leq n_0$. For this, we bound below the probability in (5.16) with the probability that, among the first $T(n, \ell)$ jumps of the renewal $\{f(k)\}_k$, there is one of length $n/2$ and all the others have length 0. In formulas,

$$\mathbb{P}(f(T(n,\ell)) \in (n/4, 3n/4)) \geq T(n,\ell)Q(0)^{T(n,\ell)-1}Q(n/2),$$

where the prefactor $T(n,\ell)$ comes from the choice of the location of the long jump. Recalling the definition of $Q(\cdot)$ and the estimate (A.5) one sees that this expression is lower bounded by a quantity which depends only on $n$ and therefore by a positive constant, since $n \leq n_0$.

## Appendix B. On the Green's function of a transient renewal

Here we prove the following statement, which is a bit more general than (6.17):

**Proposition 4.** *Let $K(\cdot)$ be a probability law on the positive integers such that $n \mapsto K(n)$ is decreasing and $K(n) \overset{n \to \infty}{\sim} \mathcal{C}_K/n^{1+\rho}$ for some $\rho \in (0,1)$ and some $\mathcal{C}_K > 0$. Let $\lambda < 1$ and $K_\lambda(\cdot) := \lambda K(\cdot)$. Then, if $\mathbf{P}_\lambda$ is the law of the transient renewal $\hat{S}$ with inter-arrival law $K_\lambda(\cdot)$, one has*

$$\mathbf{P}_\lambda(N \in \hat{S}) \overset{N \to \infty}{\sim} \frac{K_\lambda(N)}{(1-\lambda)^2}(1 + O(N^{-w})) \tag{B.1}$$

*for some $0 < w < \rho$.*

Without the explicit estimate on the error term this result is well known, see e.g. [10, Th. A.4].

*Proof.* Let $\mathbf{P} := \mathbf{P}_{\lambda=1}$. Start by writing

$$\mathbf{P}_\lambda(N \in \hat{S}) = \sum_{k=1}^{N} \lambda^k \mathbf{P}(\hat{S}_k = N) \tag{B.2}$$

and remark that if $C = C(\lambda)$ is sufficiently large, the terms with $k > C \log N$ are negligible and give a contribution which is say $o(N^{-1-2\rho}) = o(K(N)N^{-w})$. Next we show that

$$\mathbf{P}(\hat{S}_k = N) = kK(N)(1 + O(k^c/N^w)) \tag{B.3}$$

for some positive $c$ and some $0 < w < \rho$, which together with the previous observation about the large values of $k$ implies the desired bound (B.2). The lower bound is easy:

$$\mathbf{P}(\hat{S}_k = N) \geq k\mathbf{P}(\hat{S}_k = N, \hat{S}_1 > (3/4)N) \geq kK(N) \sum_{x=(3/4)N}^{N} \mathbf{P}(\hat{S}_{k-1} = N - x)$$

$$= kK(N)\{1 - \mathbf{P}(\hat{S}_{k-1} > N/4)\}$$

$$\geq kK(N)\left\{1 - (k-1)\mathbf{P}\left[\hat{S}_1 > N/(4(k-1))\right]\right\}$$

$$\geq kK(N)\left(1 - \frac{c(k-1)^{1+\rho}}{N^\rho}\right)$$

In the first inequality, we used the partition of $\{\hat{S}_k = N\}$ into the sets $\{\hat{S}_k = N\} \cap \{\hat{S}_i - \hat{S}_{i-1} > (3/4)N\}$, $i = 1, \ldots, k$, and the exchangeability of the law of $\{\hat{S}_i - \hat{S}_{i-1}\}_{i \leq k}$. The second one follows from that $K(\cdot)$ is monotone decreasing.

Let us now prove the upper bound. For any fixed $0 < \delta < 1$, by the same argument as above, we have

$$
\begin{aligned}
\mathbf{P}(\hat{S}_k = N) = {} & k\,\mathbf{P}(\{\hat{S}_k = N\} \cap \{\hat{S}_1 \geq (N - N^\delta)\}) \\
& + \mathbf{P}(\{\hat{S}_k = N\} \cap \{\nexists i \leq k \,:\, \hat{S}_i - \hat{S}_{i-1} \geq (N - N^\delta)\}).
\end{aligned}
\tag{B.4}
$$

On one hand one has

$$
\begin{aligned}
k\,\mathbf{P}[\hat{S}_k = N, \hat{S}_1 \geq (N - N^\delta)] &= k \sum_{x=(N-N^\delta)}^{N} K(x)\,\mathbf{P}[\hat{S}_{k-1} = N - x] \\
&\leq k\,K(N - N^\delta) = k\,K(N)(1 + O(N^{\delta-1}))
\end{aligned}
\tag{B.5}
$$

(there is at most one jump longer than $N - N^\delta$ if $N$ is large). On the other hand,

$$
\begin{aligned}
& \mathbf{P}(\hat{S}_k = N, \nexists i : (\hat{S}_i - \hat{S}_{i-1}) \geq (N - N^\delta))) \\
& \quad \leq k\mathbf{P}(\hat{S}_k = N; (N - N^\delta) \geq \hat{S}_1 \geq (\hat{S}_i - \hat{S}_{i-1}), i = 2, \ldots, k) \\
& \quad = k \sum_{x=(N/k)}^{N-N^\delta} \mathbf{P}(\hat{S}_k = N; x = \hat{S}_1 \geq (\hat{S}_i - \hat{S}_{i-1}), i = 2, \ldots, k) \\
& \quad \leq k \sum_{x=(N/k)}^{N-N^\delta} K(x)\mathbf{P}(\hat{S}_{k-1} = N - x) \leq kK(N/k)\mathbf{P}(\hat{S}_{k-1} \geq N^\delta) \\
& \quad \leq k(k-1)K(N/k)\mathbf{P}(\hat{S}_1 \geq N^\delta/(k-1)) \leq c\,k^c K(N) \times O(N^{-\delta\rho}).
\end{aligned}
\tag{B.6}
$$

$\square$

**Acknowledgments** We would like to thank Fabio Martinelli for enlightening conversations and for suggesting the proof of Lemma 5. Both authors acknowledge the support of Agence Nationale de la Recherche through grant ANR-2010-BLAN-0108. F.T. was supported also by European Research Council through the "Advanced Grant" PTRELSS 228032.

Université de Lyon, CNRS (UMPA)
Ecole Normale Supérieure de Lyon,
46, allée d'Italie,
69364 Lyon Cedex 07 - France.
`Cedric.Bernardin@umpa.ens-lyon.fr`
`http://www.umpa.ens-lyon.fr/~cbernard/`

Université de Lyon and CNRS
Ecole Normale Supérieure de Lyon,
Laboratoire de Physique,
46, allée d'Italie,
69364 Lyon Cedex 07 - France.
`fabio-lucio.toninelli@ens-lyon.fr`